\documentclass[a4paper,10pt]{article}
\usepackage{authblk}

\usepackage[ngerman,english]{babel}
\usepackage[latin1]{inputenc}
\usepackage{csquotes}
\usepackage{amsfonts,amsmath,amsthm}
\usepackage{empheq}
\usepackage{cancel}
\usepackage[titletoc,title]{appendix}
\usepackage[backend=bibtex8,doi=false,eprint=false,firstinits=true,isbn=false,style=numeric-comp,maxnames=99]{biblatex}
\makeatletter
\def\blx@maxline{77}
\makeatother

\bibliography{Bib_ECM_Degradation.bib}

\usepackage{cases}
\usepackage{mathabx}
\usepackage{bbm}
\usepackage{xfrac}
\usepackage{fancyhdr}
\usepackage{color}
\usepackage[colorinlistoftodos,textsize=small]{todonotes}
\usepackage[colorlinks]{hyperref}
\definecolor{blue75}{rgb}{0,0,.75}
\definecolor{green75}{rgb}{0,.75,0}
\hypersetup{colorlinks=true, urlcolor=blue75,linkcolor=blue75,citecolor=green75,pdfstartview=FitB,bookmarksopen=true,bookmarksopenlevel=1}
\usepackage[a4paper, left=2.5cm, right=2.0cm, top=2.5cm,bottom=2cm]{geometry}
\usepackage{constants}
\newcommand{\parenthezises}[1]{\arabic{#1}}
\newconstantfamily{C}{
	symbol=C,
	format=\parenthezises,
	reset={section}
}
\newconstantfamily{M}{
	symbol=M,
	format=\parenthezises,
	reset={section}
}
\newconstantfamily{B}{
	symbol=B,
	format=\parenthezises,
	reset={section}
}
\newconstantfamily{xi}{
	symbol=\xi,
	format=\parenthezises,
	reset={section}
}
\usepackage{enumerate}
\usepackage{graphicx}
\graphicspath{ {images/} }
\usepackage{wrapfig}
\usepackage{figbib}
\allowdisplaybreaks
\usepackage[capitalise]{cleveref}

\crefdefaultlabelformat{{\it #2#1#3}}

\crefname{equation}{}{}

\crefname{enumi}{}{}
\creflabelformat{enumi}{{(#2#1#3)}}

\crefname{section}{{\it Section}}{{\it Sections}}
\crefname{subsection}{{\it Subsection}}{{\it Subsections}}
\crefname{subsubsection}{{\it Paragraph}}{{\it Paragraphs}}
\newtheorem{Theorem}{Theorem}[section]
\crefname{Theorem}{{\it Theorem}}{{\it Theorems}}
\newtheorem{Definition}[Theorem]{Definition}
\crefname{Definition}{{\it Definition}}{{\it Definitions}}
\newtheorem{Lemma}[Theorem]{Lemma}
\crefname{Lemma}{{\it Lemma}}{{\it Lemmas}}

\crefname{Proposition}{{\it Proposition}}{{\it Propositions}}

\crefname{Assumption}{{\it Assumption}}{{\it Assumptions}}

\crefname{Corollary}{{\it Corollary}}{{\it Corollaries}}

\theoremstyle{definition}

\crefname{Remark}{{\it Remark}}{{\it Remarks}}

\crefname{Notation}{{\it Notation}}{{\it Notations}}

\crefname{Example}{{\it Example}}{{\it Examples}}

\newcommand{\be}{\begin{equation} \label}
\newcommand{\ee}{\end{equation}}
\newcommand{\bea}{\begin{eqnarray}\label}
\newcommand{\eea}{\end{eqnarray}}
\newcommand{\bas}{\begin{eqnarray*}}
	\newcommand{\eas}{\end{eqnarray*}}
\newcommand{\bit}{\begin{itemize}}
	\newcommand{\eit}{\end{itemize}}
\newcommand{\R}{\mathbb{R}}

\newcommand{\cb}{\color{blue}}
\begin{document}
	\enlargethispage{10mm}

\title{Global existence of weak solutions to a cell migration and (de)differentiation model with double haptotaxis in the context of tissue regeneration}
\author{Nishith Mohan\footnote{mohan@mathematik.uni-kl.de}\quad  and\quad  Christina Surulescu\footnote{surulescu@mathematik.uni-kl.de}\\
	{\small RPTU Kaiserslautern-Landau, Department of Mathematics,} \\
	{\small Gottlieb-Daimler-Str. 48, 67663 Kaiserslautern, Germany}}

\date{\today}
\maketitle

\begin{abstract} We study a model for the spread and (de)differentiation of mesenchymal stem cells and chondrocytes in a scaffold whose fibers are coated with hyaluron. The chondrocytes produce new extracellular matrix, which, together with hyaluron, serves as haptotactic cue for the stem cell migration. We prove global existence of weak solutions of the corresponding cross-diffusion system with double haptotaxis. 
\end{abstract}	

\section{Introduction and model set-up}

Tissue regeneration is a continuously expanding research field, with the ultimate aim of identifying the essential factors involved in cartilage or bone repair after traumas or degeneration of joints, and of conceiving best suitable treatment approaches. The interplay between mathematical modeling, biology, medicine, and biomedical engineering holds the promise for reaching this aim. We refer to \cite{Waters2021} for a review on such multidisciplinary works and to \cite{Grosjean} for a review of mathematical model classes employed in this context.\\[-2ex]

\noindent
In this note we consider a simplified version of the model derived in \cite{Grosjean} by a multiscale approach and of the model in \cite{paper-pamm}, where data-based numerical simulations have been performed, showing biologically reasonable behavior triggered by the same kind of haptotaxis as the one considered here. 
Concretely, we investigate the global existence of solutions to a system characterizing the dynamics of mesenchymal stem cells (MSCs) seeded in a scaffold whose fibres are coated with hyaluron. They differentiate into chondrocytes, who are able to express collagen containing extracellular matrix (ECM), thus rebuilding tissue; they can dedifferentiate back into stem cells. MSCs are attracted towards gradients of hyaluron and ECM. As none of this tactic cues is diffusing, MSCs perform, beyond diffusion,  double haptotaxis. Their proliferation is limited in a logistic manner by the densities of MSC, chondrocytes, and ECM. Chondrocytes can diffuse, arise from differentiated MSCs, and can de-differentiate. Hyaluron can be uptaken by MSCs and chondrocytes, infer a natural decay, and be produced by chiondrocytes to a limited amount. ECM can be degraded by MSCs (e.g., by expression of matrix degrading proteins), infer natural decay and (limited) production by chondrocytes. \\[-2ex]

\noindent
Let us denote by $c_1$ and $c_2$ the densities of mesenchymal stem cells and chondrocytes, respectively, with $h$ the concentration of hyaluron impregnating scaffold fibers, and with $\tau $ the density of newly produced ECM. All these quantities are space- and time-dependent. With these variables we study the following system:
\begin{equation}
  \begin{cases}
  \label{nondim1}
    \partial_t c_1 = a_1\Delta c_{1} - \nabla \cdot (b_h c_1 \nabla h) - \nabla \cdot (b_\tau c_1 \nabla \tau) - \alpha_1(\tau)c_1 + \alpha_2(\tau)c_2 \\[7pt]
    \hspace{8cm}+ \beta c_{1}(1 - c_{1} - c_{2} - \tau),                                    &  x \in \Omega,\ t > 0,           \\[7pt]
    \partial_t c_2  = a_2 \Delta c_2 + \alpha_1(\tau)c_1 -\alpha_2( \tau)c_2,                   &  x \in \Omega,\ t > 0,           \\[7pt]
    \partial_t h = - \gamma_1 h c_1 - \gamma_2 h c_2 - \mu h +  \dfrac{c_2}{1 + c_2},                                &  x \in \Omega,\ t > 0,           \\[10pt]
    \partial_t \tau = - \delta c_1 \tau - \sigma\tau +  \dfrac{c_2}{1 + c_2},                                             &  x \in \Omega,\ t > 0,
  \end{cases}
\end{equation}
subject to the following zero-flux boundary conditions ($\nu$ denotes the outward unit normal on the boundary of $\Omega$ )
\begin{equation}
\label{bco}
-a_1\dfrac{\partial c_1}{\partial \nu} + b_h c_1 \dfrac{\partial h}{\partial \nu} + b_\tau c_1 \dfrac{\partial \tau}{\partial \nu} = \dfrac{\partial c_2}{\partial \nu} = 0, \quad x \in \partial \Omega,\ t > 0,
\end{equation}
and the initial conditions
\begin{equation}
  \label{oic}
  c_1(x, 0) = c_{10}(x) , ~~ c_2(x, 0) = c_{20}(x), ~~ h(x, 0) = h_0(x), ~~ \tau(x, 0) = \tau_0(x), \quad x \in \Omega.
\end{equation}

%%%%%%%%%%%%%%%%%%%
\noindent
The above system features double haptotaxis in the motility of MSCs. As such, it belongs to category (ii) (with missing chemotaxis) in the classification given in the review \cite{Kolbe2021}. Haptotaxis models are known to have a more challenging analysis, due to the lack of regularity caused by the degenerated diffusion of the tactic cues. Most models with multiple haptotaxis distribute the corresponding advection terms among the equations for the solution components, i.e. populations of cells \cite{Sfakianakis2016} or virions (see \cite{Alzahrani2019} and subsequent versions and analyses thereof, e.g.,   \cite{Daba2024,tao2021global,Ren2021,Tao2014,Wang2022,Zhou2025}). In \cite{Kolbe2021} a model with a douby haptotactic population of motile cells was proposed, which performed haptotaxis towards ECM gradients and away from gradients of resting and proliferating cells, thus with unilateral interspecific repellence. Global existence of weak solutions was proved by constructing a sequence of regularized approximations featuring diffusion in all equations, showing global existence of solutions thereof, deducing estimates stemming from an entropy-type functional, and then passing to the limits to show global existence of solutions to the original problem. The proof of global existence for our model \eqref{nondim1}, \eqref{bco}, \eqref{oic} with double attractive haptotaxis follows the same scheme.\\[-2ex]

%%%%%%%%%%%%%%%%%%%%%%%%%%	

\noindent
In the sequel we assume that  $n \in \{2,3\}$ and $\Omega \subset \mathbb{R}^n$  is a bounded convex domain with a smooth enough  boundary. Moreover, we assume that
\begin{equation}
\begin{cases}
  \label{oic2}
  c_{10},  c_{20} \in C^0(\bar{\Omega}), ~~h_{0}, \tau_{0} \in W^{1, 2}(\Omega) \cap  C^0(\bar{\Omega}),  ~ \sqrt{h_0},\  \sqrt{\tau_0} \in W^{1, 2}(\Omega),  \\[5pt]
  c_{10},\ c_{20} \geq 0,\ h_0,\ \tau_0 > 0 ~\text{in}~ \Omega ~\text{with}~ c_{10} \not\equiv 0 ~\text{and}~ c_{20} \not\equiv 0,
\end{cases}
\end{equation}
and the functions $\alpha_i$, $i \in \{1, 2\}$ satisfy
\begin{equation}
\label{alpha_def}
\begin{cases}
  & \alpha_{i}(z) \in C^{\vartheta, \frac{\vartheta}{2}}(\bar{\Omega} \times [0, T]),\quad (\vartheta \in (0, 1) ~\text{and}~ T > 0), \quad \alpha_i(z) > 0, ~\text{and} \\[5pt]
  & \alpha_{i}(z) \leq M_{\alpha_i}, \quad i = \{1, 2\},  ~\text{for all}~ z \geq 0
\end{cases}
\end{equation}
and all the parameters $a_1, a_2, b_h, b_\tau, \beta, \gamma_1, \gamma_2, \delta, \mu$ and $\sigma$ are positive.\\[-2ex]

\noindent
The primary objective of this work is to prove that problem \eqref{nondim1}-\eqref{oic} admits global weak solutions, in the sense given by the following definition:
\begin{Definition}
  \label{solution_definition}
  Let $T \in (0, \infty)$. A weak solution to the problem \eqref{nondim1}-\eqref{oic} in $\Omega \times (0, T)$ consists of a quadruple of nonnegative functions $(c_1, c_2, h, \tau)$ such that
  \begin{equation*}
    \begin{cases}
      & c_1 \in L^2(\Omega \times (0, T)) \cap L^{\frac{4}{3}}((0, T); W^{1, \frac{4}{3}}(\Omega)),             \\[5pt]
      & c_2 \in L^{\frac{5}{4}}((0, T); W^{1, \frac{5}{4}}(\Omega)),             \\[5pt]
      & h \in L^\infty(\Omega \times (0, T)) \cap L^2((0, T); W^{1, 2}(\Omega)),   \quad ~\text{and}     \\[5pt]
      & \tau \in L^\infty(\Omega \times (0, T)) \cap L^2((0, T); W^{1, 2}(\Omega)),
    \end{cases}
  \end{equation*}
  and the folowing equations are satisfied:
\begin{align}
  \label{defeq1}
  \nonumber  & - \int_0^T \int_\Omega c_{1} \partial_t \phi - \int_\Omega c_{10} \phi(\cdot, 0) = - a_1 \int_0^T \int_\Omega \nabla c_{1} \cdot \nabla \phi + b_h \int_0^T \int_\Omega c_{1} \nabla h \cdot \nabla \phi \\[5pt]
  \nonumber & + b_\tau \int_0^T \int_\Omega c_{1} \nabla \tau \cdot \nabla \phi - \int_0^T \int_\Omega \alpha_1(\tau)c_{1} \phi + \int_0^T \int_\Omega \alpha_2(\tau)c_{2} \phi  \\[5pt]
  & + \beta \int_0^T \int_\Omega c_{1}(1 - c_{1}  -  c_{2}  - \tau ) \phi
\end{align}
and
\begin{align}
  \label{defeq2}
- \int_0^T \int_\Omega c_{2} \partial_t \phi - \int_\Omega c_{20} \phi(\cdot, 0) = &- a_2 \int_0^T \int_\Omega \nabla c_{2} \cdot \nabla \phi + \int_0^T \int_\Omega \alpha_1(\tau)c_{1} \phi - \int_0^T \int_\Omega \alpha_2(\tau) c_{2} \phi
\end{align}
and
\begin{align}
  \label{defeq3}
 - \int_0^T \int_\Omega h \partial_t \phi - \int_\Omega h_{0} \phi(\cdot, 0) & = - \gamma_1 \int_0^T \int_\Omega h c_{1} \phi  - \gamma_2 \int_0^T \int_\Omega h c_{2}\phi - \mu \int_0^T \int_\Omega  h \phi  +  \int_0^T \int_\Omega \dfrac{c_{2}}{1 + c_{2}} \phi ,
\end{align}
as well as
\begin{align}
  \label{defeq4}
  - \int_0^T \int_\Omega \tau \partial_t \phi - \int_\Omega \tau_{0} \phi(\cdot, 0) = - \delta \int_0^T \int_\Omega \tau c_{1} \phi - \sigma \int_0^T \int_\Omega \tau \phi - \int_0^T \int_\Omega \dfrac{c_{2}}{1 + c_{2}} \phi
\end{align}
for all $\phi \in C^\infty_0(\bar{\Omega} \times [0, T))$ with $\frac{\partial \phi}{\partial \nu} = 0$ on $\partial \Omega \times [0, T)$. If the quadruple $(c_1, c_2, h, \tau)$ is a weak solution to \eqref{nondim1}-\eqref{oic} in $\Omega \times (0, T)$ for all $T > 0$, then it is referred to as a global weak solution.
\end{Definition}

\noindent
Our main result asserts that the problem \eqref{nondim1}-\eqref{oic} admits a global weak solution:

\begin{Theorem}
  \label{main_theorem}
Assume that \eqref{oic2} holds true and the transition functions $\alpha_i$, where $i \in \{1, 2\}$, satisfy \eqref{alpha_def}. Then, the problem \eqref{nondim1}-\eqref{oic} has at least one global weak solution in the sense of Definition \ref{solution_definition}.
\end{Theorem}

\noindent
In the sequel we will use the following notations and conventions:
\begin{itemize}
	\item The integrals $\int_\Omega f(x) dx$ are abbreviated as $\int_\Omega f(x)$.
	\item The sequentiality of the constants $C_i, i = 1, 2, 3, \ldots$ holds only within the lemma/theorem and its proof in which the constants are used. The sequence restarts once the proof is over.
	\item \eqref{alpha_def} holds true throughout the discussion.
\end{itemize}

\section{A family of approximate problems}
In order to construct a solution for \eqref{nondim1}, \eqref{bco} and \eqref{oic}, we first define a family $\{F_\varepsilon\}_{\varepsilon \in (0, 1)} \subset C^\infty([0, \infty))$ of functions such that for any $\varepsilon \in (0, 1)$, we have
\begin{equation}
  \label{app1}
  F_\varepsilon(s) :=  \frac{s}{1 + \varepsilon s} \quad \text{for} ~ \varepsilon \in (0, 1) ~\text{and}~ s \geq 0.
\end{equation}
Clearly, for every $\varepsilon \in (0, 1)$ we have
\begin{align}
  \label{app2} & 0 \leq F_\varepsilon(s) \leq s,                                   \quad \text{for all} ~ s \geq 0, \quad \text{and} \\[5pt]
  \label{app3} & 0 \leq F_\varepsilon(s) < \frac{1}{\varepsilon}                   \quad \text{for all} ~ s \geq 0.
\end{align}
We will consider the following regularized problems to construct a weak solution for \eqref{nondim1}-\eqref{oic}:
\begin{equation}
  \begin{cases}
  \label{model1}
    \partial_t c_{1\varepsilon } = a_1 \Delta c_{1\varepsilon} - \nabla \cdot (b_h c_{1 \varepsilon} \nabla h_\varepsilon) - \nabla \cdot (b_\tau c_{1\varepsilon} \nabla \tau_\varepsilon) - \alpha_1(\tau_\varepsilon)c_{1\varepsilon} + \alpha_2( \tau_\varepsilon) F_\varepsilon(c_{2\varepsilon}) \\[5pt]
    \hspace{3cm}+ \beta c_{1\varepsilon}(1 - c_{1\varepsilon} -  c_{2\varepsilon} - \tau_\varepsilon)  - \varepsilon c_{1\varepsilon}^\theta,                                & \hspace{-1cm} x \in \Omega, t > 0,           \\[5pt]
    \partial_t c_{2 \varepsilon} = a_2\Delta c_{2\varepsilon} + \alpha_1(\tau_\varepsilon) c_{1\varepsilon} -  \alpha_2 (\tau_\varepsilon) F_\varepsilon(c_{2\varepsilon}),           & \hspace{-1cm} x \in \Omega, t > 0,           \\[5pt]
    \partial_t h_\varepsilon  = \varepsilon \Delta h_\varepsilon - \gamma_1 h_\varepsilon c_{1\varepsilon}  - \gamma_2 h_\varepsilon c_{2\varepsilon} - \mu h_\varepsilon +  \dfrac{c_{2\varepsilon}}{1 + c_{2\varepsilon}}  & \hspace{-1cm}  x \in \Omega, t > 0, \\[10pt]
    \partial_t \tau_\varepsilon = \varepsilon \Delta \tau_\varepsilon - \delta \tau_\varepsilon c_{1\varepsilon}  - \sigma \tau_\varepsilon +  \dfrac{c_{2\varepsilon}}{1 + c_{2\varepsilon}} ,   &\hspace{-1cm} x \in \Omega, t > 0,    \\[7pt]
    \partial_\nu c_{1\varepsilon} = \partial_\nu c_{2\varepsilon} = \partial_\nu h_\varepsilon = \partial_\nu \tau_\varepsilon = 0,                                             &\hspace{-1cm}  x \in \partial \Omega, t > 0 , \\[5pt]
    c_{1\varepsilon}(x, 0) = c_{10\varepsilon}(x), c_{2\varepsilon}(x, 0) = c_{20\varepsilon}(x), h_{\varepsilon}(x, 0) = h_{0\varepsilon}(x), \tau_{\varepsilon}(x, 0) = \tau_{0\varepsilon}(x),                     & \hspace{-1cm}  x \in \Omega,
  \end{cases}
\end{equation}
for $\varepsilon \in (0, 1)$, where $\theta > \max\{2, n\}$ is a fixed parameter and $F_\varepsilon(c_{i\varepsilon}), i \in \{1, 2\}$ is as defined in \eqref{app1} satisfying \eqref{app2} and \eqref{app3}. The families of functions $\{c_{10\varepsilon}\}_{\varepsilon \in (0, 1)}, \{c_{20\varepsilon}\}_{\varepsilon \in (0, 1)}, \{h_{0\varepsilon}\}_{\varepsilon \in (0, 1)}$ and $\{\tau_{0\varepsilon}\}_{\varepsilon \in (0, 1)}$ satisfy
\begin{equation}
\label{reg_assumptions}
  \begin{cases}
    c_{10\varepsilon}, c_{20\varepsilon}, h_{0\varepsilon}, \tau_{0\varepsilon} \in C^3(\bar{\Omega}),                                                                   \\[3pt]
    c_{10\varepsilon} > 0, c_{20\varepsilon} > 0, h_{0\varepsilon} > 0, \tau_{0\varepsilon} > 0 ~\text{in}~ \bar{\Omega},                                                \\[3pt]
    \partial_\nu c_{10\varepsilon} = \partial_\nu c_{20\varepsilon} = \partial_\nu h_{0\varepsilon} = \partial_\nu \tau_{0\varepsilon} = 0 ~\text{on}~ \partial \Omega,  \\[3pt]
    c_{10\varepsilon} \to c_{10}, c_{20\varepsilon} \to c_{20},  ~\text{in}~ C^0(\bar{\Omega}) ~\text{as}~ \varepsilon \searrow 0,                       \\[3pt]
    \sqrt{h_{0\varepsilon}} \to \sqrt{h_0} ~\text{and}~ \sqrt{\tau_{0\varepsilon}} \to \sqrt{\tau_0} ~\text{in}~ W^{1, 2}(\Omega) \cap C^0(\bar{\Omega}),  ~\text{as}~ \varepsilon \searrow 0.
  \end{cases}
\end{equation}

\section{Global existence for approximate problems}

\begin{Lemma}
\label{lemma_loc_exist}
Assume that \eqref{reg_assumptions} holds true. Then, for every $\varepsilon \in (0, 1)$ there exist $T_{\max, \varepsilon} \in (0, \infty]$ and a collection of positive functions $c_{1\varepsilon}, c_{2\varepsilon}, h_\varepsilon$ and $\tau_\varepsilon$, each belonging to $C^{2,1}(\bar{\Omega} \times [0, T_{\max, \varepsilon}))$ such that the quadruple $(c_{1\varepsilon}, c_{2\varepsilon}, h_\varepsilon, \tau_\varepsilon)$ solves \eqref{model1} classically in $\Omega \times (0, T_{\max, \varepsilon})$. Moreover, if $~ T_{\max, \varepsilon} < \infty$, then
  \begin{equation}
  \label{loe1}
    \limsup_{t \nearrow T_{\max, \varepsilon}} \left\{\|c_{1\varepsilon}(\cdot, t)\|_{C^{2 + \vartheta}(\bar{\Omega})} + \|c_{2\varepsilon}(\cdot, t)\|_{C^{2 + \vartheta}(\bar{\Omega})} + \|h_{\varepsilon}(\cdot, t)\|_{C^{2 + \vartheta}(\bar{\Omega})} + \|\tau_{\varepsilon}(\cdot, t)\|_{C^{2 + \vartheta}(\bar{\Omega})}\right\} = \infty
  \end{equation}
 for all $\vartheta \in (0, 1)$.
\end{Lemma}
\begin{proof}
By adapting the arguments from \cite[Lemma 3.1]{stinner2014global}, we can easily prove this result.
\end{proof}

%%%%%%%%%%%%%%%%%%%%%%%%%%%%%%%%%%%%%%%%%%%%%%%%%%%%%%%%%%%%%%%%%%%%%%%%%%%%%%%%%%%%%%%%%%%%%%%%%%%%%%%%%%%%%%%%%%%%%%%%%%%%%%%%%%%%%%%%%%%%%%%%%%%%%%%%%%%%%%%%%%%%%%%%%%%%%%%%%%%%%%%%%%%%%%%%%%%%%%%%%%%%%%%

\begin{Lemma}
  \label{lemmale1}
  Let $T > 0$. Then there exists a $C(T) > 0$ such that for all $\varepsilon \in (0, 1)$
  \begin{equation}
    \label{L1bound}
    \int_\Omega c_{1\varepsilon}(\cdot, t) \leq C(T) ~~\text{and}~~  \int_\Omega c_{2\varepsilon}(\cdot, t) \leq C(T),
  \end{equation}
  for all $t \in (0, \widehat{T}_\varepsilon)$. Additionally, we have
    \begin{equation}
    \label{sqL1bound}
    \int_0^{\widehat{T}_\varepsilon} \int_\Omega c_{1\varepsilon}^2 \leq C(T),
  \end{equation}
  \begin{equation}
    \label{epL1bound}
     \varepsilon \int_0^{\widehat{T}_\varepsilon} \int_\Omega c_{1\varepsilon}^\theta \leq C(T),
  \end{equation}
   where $\widehat{T}_\varepsilon := \min\{T, T_{\max, \varepsilon}\}$.
\end{Lemma}

\begin{proof}
  We integrate with respect to space the first and second equation in \eqref{model1}, then add the respective results. This yields
  \begin{align}
    \label{L1bound1}
     \dfrac{d}{dt} \left\{\int_\Omega c_{1\varepsilon} + \int_\Omega c_{2\varepsilon}\right\} & = \beta \int_\Omega c_{1\varepsilon} (1 -  c_{2\varepsilon} - \tau_\varepsilon) - \beta \int_\Omega c_{1\varepsilon}^2  - \varepsilon  \int_\Omega c_{1\varepsilon}^\theta \quad \text{for all} ~t \in (0,  T_{\max, \varepsilon})
    \end{align}
 As the solution components are non-negative, we can have
 \begin{align}
   \label{L1bound2}
   \dfrac{d}{dt} \left\{\int_\Omega c_{1\varepsilon} + \int_\Omega c_{2\varepsilon} \right\}  + \beta \int_\Omega c_{1\varepsilon}^2 +\varepsilon  \int_\Omega c_{1\varepsilon}^\theta& \leq \beta  \left\{ \int_\Omega c_{1\varepsilon} + \int_\Omega c_{2\varepsilon} \right\}  \quad \text{for all} ~ t \in (0, T_{\max, \varepsilon}).
 \end{align}
 Clearly,
 \begin{align}
   \label{L1bound3}
   \dfrac{d}{dt} \left\{\int_\Omega c_{1\varepsilon} + \int_\Omega c_{2\varepsilon} \right\} & \leq \beta  \left\{ \int_\Omega c_{1\varepsilon} + \int_\Omega c_{2\varepsilon} \right\} \quad \text{for all} ~ t \in (0, T_{\max, \varepsilon}).
 \end{align}
 Applying Gronwall's inequality to \eqref{L1bound3} over the interval $t \in (0, \widehat{T}_\varepsilon)$ leads to
 \begin{equation}
   \label{L1bound4}
  \int_\Omega c_{1\varepsilon}(\cdot, t) + \int_\Omega c_{2\varepsilon}(\cdot, t)  \leq e^{C_1 \widehat{T}_\varepsilon} \left\{\sup_{\varepsilon \in (0, 1)}\int_\Omega\left(c_{10\varepsilon} + c_{20\varepsilon}\right) + C_2 \widehat{T}_\varepsilon \right\}.
 \end{equation}
This establishes \eqref{L1bound}. In light of \eqref{L1bound}, we can obtain \eqref{sqL1bound} and \eqref{epL1bound} by integrating \eqref{L1bound2} over $(0, \widehat{T}_\varepsilon)$.

\end{proof}

%%%%%%%%%%%%%%%%%%%%%%%%%%%%%%%%%%%%%%%%%%%%%%%%%%%%%%%%%%%%%%%%%%%%%%%%%%%%%%%%%%%%%%%%%%%%%%%%%%%%%%%%%%%%%%%%%%%%%%%%%%%%%%%%%%%%%%%%%%%%%%%%%%%%%%%%%%%%%%%%%%%%%%%%%%%%%%%%%%%%%%%%%%%%%%%%%%%%%%%%%%%%%%%

\begin{Lemma}
  \label{ht_bound}
   For all $\varepsilon \in (0, 1)$ the third and fourth solution components of \eqref{model1} satisfy
  \begin{align}
    \label{h_bound} & \|h_\varepsilon(\cdot, t)\|_{L^\infty(\Omega)} \leq \frac{r_\ast}{\mu} + \|h_{0 \varepsilon }\|_{L^\infty(\Omega)} =: M_h      \quad \text{for all}~ t \in (0, T_{\max, \varepsilon}) \\[5pt]
    \label{t_bound} & \|\tau_\varepsilon(\cdot, t)\|_{L^\infty(\Omega)} \leq \frac{r_\ast}{\sigma} + \|\tau_{0 \varepsilon }\|_{L^\infty(\Omega)} =: M_\tau \quad \text{for all}~ t \in (0, T_{\max, \varepsilon}),
  \end{align}
   where $r:= \dfrac{c_{2\varepsilon}}{1 + c_{2\varepsilon}}$ and $r_\ast := \|r\|_{L^\infty(\Omega \times (0, \infty))}$.
\end{Lemma}

\begin{proof}
  Let
  \begin{equation*}
    \bar{h}(t) := \|h_{0\varepsilon}\|_{L^\infty(\Omega)} e^{- \mu t} + \int_0^t e^{- \mu (t - s)} \|r(\cdot, s)\|_{L^\infty(\Omega)} ds
  \end{equation*}
  for $t \in (0, \infty)$. Then %{\cb {\cancel{$\bar{h}$ is spatially homogenous for each $t \in (0, T_{\max, \varepsilon})$ with}}} 
  $\bar{h} \leq M_h$ on $\bar{\Omega} \times (0, T_{\max, \varepsilon})$ and it holds
  \begin{align*}
    \bar{h}_t - \varepsilon \Delta \bar{h} + \gamma_1 c_{1\varepsilon} \bar h + \gamma_2 c_{2\varepsilon} \bar h + \mu \bar{h}- r(x, t) \geq \bar{h}_t + \mu \bar{h} - \|r(\cdot, t)\|_{L^\infty(\Omega)}\ge  0
  \end{align*}
  on  $\bar{\Omega} \times (0, T_{\max, \varepsilon})$, owing to the non-negativity of the solution components $c_{1\varepsilon}, c_{2\varepsilon}$, and $h_\varepsilon$. Applying an ODE comparison principle results in $\bar{h} \geq h_\varepsilon$ on $\Omega \times (0, T_{\max, \varepsilon})$ and, in particular $\|h_\varepsilon(\cdot, t)\|_{L^\infty(\Omega)} \leq M_h$ for all $t \in (0, T_{\max, \varepsilon})$. This gives \eqref{h_bound}. We can derive \eqref{t_bound} analogously.

\end{proof}
%%%%%%%%%%%%%%%%%%%%%%%%%%%%%%%%%%%%%%%%%%%%%%%%%%%%%%%%%%%%%%%%%%%%%%%%%%%%%%%%%%%%%%%%%%%%%%%%%%%%%%%%%%%%%%%%%%%%%%%%%%%%%%%%%%%%%%%%%%%%%%%%%%%%%%%%%%%%%%%%%%%%%%%%%%%%%%%%%%%%%%%%%%%%%%%%%%%%%%%%%%%%%%%
\noindent
With these estimates in hand we will  now proceed to show that for any fixed $\varepsilon \in (0, 1)$ the solution quadruple $(c_{1\varepsilon}, c_{2\varepsilon}, h_\varepsilon, \tau_\varepsilon)$ obtained in Lemma \ref{lemma_loc_exist} is global in time.

\begin{Lemma}[Global existence for \eqref{model1}]
  \label{global_existence}
  For all $\varepsilon \in (0, 1)$, the solution $(c_{1\varepsilon}, c_{2\varepsilon}, h_\varepsilon, \tau_\varepsilon)$ of \eqref{model1} obtained in Lemma \ref{lemma_loc_exist} is global, i.e. $T_{\max, \varepsilon} = \infty$.
\end{Lemma}

\begin{proof}
  We will show $T_{\max, \varepsilon} = \infty$ by contradiction. Fix $\varepsilon \in (0, 1)$ and suppose that $T_{\max, \varepsilon} < \infty$.
 From \eqref{epL1bound} we can have a $C_1(\varepsilon, T) > 0$ such that
  \begin{equation}
  \label{neq1}
    \int_0^{\widehat{T}_\varepsilon} \int_\Omega c_{1\varepsilon}^\theta \leq C_1(\varepsilon, T)
  \end{equation}
  where $\widehat{T}_\varepsilon := \min\{T, T_{\max, \varepsilon}\}$ for $T > 0$. \eqref{alpha_def}, \eqref{app3} and \eqref{neq1} allow us to apply the well known maximal Sobolev regularity for parabolic equations \cite{matthias1997heat} to the second equation  of \eqref{model1}, providing a $C_2(\varepsilon, T) > 0$ such that
  \begin{align*}
    \int_0^{\widehat{T}_\varepsilon} \|\partial_tc_{2\varepsilon }\|^\theta_{L^\theta(\Omega)} + \int_0^{\widehat{T}_\varepsilon} \|c_{2\varepsilon }\|^\theta_{W^{2, \theta}(\Omega)} & \leq C_2\|c_{20\varepsilon}\|^\theta_{W^{2, \theta}(\Omega)} + C_2 \int_0^{\widehat{T}_\varepsilon} \left (\|\alpha_1(\tau_\varepsilon) c_{1\varepsilon}\|^\theta_{L^\theta(\Omega)} + \|\alpha_2(\tau_\varepsilon) F_\varepsilon(c_{2\varepsilon})\|^\theta_{L^\theta(\Omega)}\right )  \\[5pt]
    & \leq C_2\left(\|c_{20\varepsilon}\|^\theta_{W^{2, \theta}(\Omega)} + M_{\alpha_1}^\theta \int_0^{\widehat{T}_\varepsilon} \|c_{1\varepsilon}\|^\theta_{W^{2, \theta}(\Omega)} + \varepsilon^{-\theta}M_{\alpha_2}^\theta\widehat{T}_\varepsilon \right) \\[5pt]
    & \leq C_2(\varepsilon, T).
  \end{align*}
  From above we can conclude that
  \begin{align}
    \label{sob_c2}
    \int_0^{\widehat{T}_\varepsilon} \|c_{2\varepsilon }\|^\theta_{W^{2, \theta}(\Omega)} \leq C_2(\varepsilon, T).
  \end{align}
Using \eqref{h_bound}, \eqref{neq1}, \eqref{sob_c2}, and the fact that $\frac{c_{2\varepsilon}}{1 + c_{2\varepsilon}} \leq 1$ for $c_{2\varepsilon} \geq 0$, we can show that $f_\varepsilon := \partial_t h_\varepsilon - \varepsilon \Delta h_\varepsilon \equiv - \gamma_1 h_\varepsilon c_{1\varepsilon}  - \gamma_2 h_\varepsilon c_{2\varepsilon} - \mu h_\varepsilon +  \tfrac{c_{2\varepsilon}}{1 + c_{2\varepsilon}}$ is bounded in $L^\theta(\Omega \times (0, \widehat{T}_\varepsilon))$. Again, by maximal Sobolev regularity applied to the third equation of \eqref{model1} we can have a $C_3(\varepsilon, T) > 0$ such that
  \begin{equation}
  \label{ge1}
    \int_0^{\widehat{T}_\varepsilon} \|h_\varepsilon(\cdot, t)\|^\theta_{W^{2, \theta}(\Omega)} dt \leq C_3(\varepsilon, T).
  \end{equation}
 Similarly, $g_\varepsilon := \partial_t \tau_\varepsilon - \varepsilon \Delta \tau_\varepsilon \equiv - \delta \tau_\varepsilon c_{1\varepsilon} - \sigma \tau_\varepsilon + \frac{c_{2\varepsilon}}{1 + c_{2\varepsilon}}$ also belongs to $L^\theta(\Omega \times (0, \widehat{T}_\varepsilon))$, thanks to \eqref{neq1}, \eqref{t_bound} and the fact that $\frac{c_{2\varepsilon}}{1 + c_{2\varepsilon}} \leq 1$ for $c_{2\varepsilon} \geq 0$. Applying the same  regularity to the fourth equation in \eqref{model1}, yields a $C_4(\varepsilon, T) > 0$ fulfilling
  \begin{equation}
  \label{ge2}
    \int_0^{\widehat{T}_\varepsilon} \|\tau_\varepsilon(\cdot, t)\|^\theta_{W^{2, \theta}(\Omega)} dt \leq C_4(\varepsilon, T).
  \end{equation}
  Since $\theta > \max\{2, n\}$, applying the Sobolev embedding $W^{2, \theta}(\Omega) \hookrightarrow W^{1, \infty}(\Omega)$ in conjunction with H\"{o}lder's inequality gives us a $\widehat{C}(\varepsilon, T) > 0$ such that
  \begin{equation}
  \label{ge3}
    \int_0^{\widehat{T}_\varepsilon} \|\nabla h_\varepsilon(\cdot, t)\|^2_{L^\infty(\Omega)} dt \leq \widehat{C}(\varepsilon, T), \quad \text{and}~ \int_0^{\widehat{T}_\varepsilon} \|\nabla \tau_\varepsilon(\cdot, t)\|^2_{L^\infty(\Omega)} dt \leq \widehat{C}(\varepsilon, T).
  \end{equation}
Using the non-negativity of solution components, we test the first equation of \eqref{model1} with $pc_{1\varepsilon}^{p - 1}$, ($p > 1$)
 \begin{align}
 \label{ge6}
   \nonumber \dfrac{d}{dt}\int_\Omega c_{1\varepsilon}^{p} & = - p(p - 1) a_1 \int_\Omega c_{1 \varepsilon}^{p - 2}|\nabla c_{1 \varepsilon}|^2 + b_h p(p - 1) \int_\Omega c_{1 \varepsilon}^{p - 1} \nabla c_{1 \varepsilon} \cdot \nabla h_\varepsilon \\ & + b_\tau p(p - 1) \int_\Omega c_{1 \varepsilon}^{p - 1} \nabla c_{1 \varepsilon} \cdot \nabla \tau_\varepsilon
   - \nonumber p  \int_\Omega \alpha_1(\tau_\varepsilon)  c_{1 \varepsilon}^{p} + p \int_\Omega \alpha_2(\tau_{\varepsilon}) F_\varepsilon(c_{2 \varepsilon}) c_{1 \varepsilon}^{p - 1} \\ & + p \beta \int_\Omega c_{1\varepsilon}^p(1 - c_{1\varepsilon} - c_{2\varepsilon} - \tau_\varepsilon) - p \varepsilon \int_\Omega c_{1\varepsilon}^{\theta + p - 1}
 \end{align}
for all $t \in (0, T_{\max, \varepsilon})$. Using Young's inequality, we will first handle the two taxis terms on the right-hand side of (\ref{ge6}):
 \begin{align}
    \label{ge7} & b_h p(p - 1) \int_\Omega c_{1 \varepsilon}^{p - 1} \nabla c_{1 \varepsilon} \cdot \nabla h_\varepsilon \leq \dfrac{p(p - 1)a_1}{4} \int_ \Omega c_{1 \varepsilon}^{p - 2} |\nabla c_{1 \varepsilon}|^2 + \dfrac{p(p - 1)}{a_1} b_h^2\int_\Omega c_{1 \varepsilon}^{p} |\nabla h_{\varepsilon}|^2,   \\
    \label{ge8} & b_\tau p(p - 1) \int_\Omega c_{1 \varepsilon}^{p - 1} \nabla c_{1 \varepsilon} \cdot \nabla \tau_\varepsilon \leq \dfrac{p(p - 1){\cb a_1}}{4} \int_ \Omega c_{1 \varepsilon}^{p - 2} |\nabla c_{1 \varepsilon}|^2 + \dfrac{p(p - 1)}{a_1} b_\tau^2\int_\Omega c_{1 \varepsilon}^{p} |\nabla \tau_{\varepsilon}|^2.
 \end{align}
 Now, we consider the term $p \int_\Omega \alpha_2(\tau_{\varepsilon}) F_\varepsilon(c_{2 \varepsilon}) c_{1 \varepsilon}^{p - 1} $. Using \eqref{alpha_def}, \eqref{app3} and Young's inequality we get
 \begin{align}
 \label{ge9}
   \nonumber p \int_\Omega \alpha_2(\tau_{\varepsilon}) F_\varepsilon(c_{2 \varepsilon}) c_{1 \varepsilon}^{p - 1} & \leq p M_{\alpha_2} \varepsilon^{-1} \int_\Omega c_{1 \varepsilon}^{p - 1} \\[5pt]
   & \leq p \int_\Omega  c_{1 \varepsilon}^{p} + C(p)  M_{\alpha_2}^p \varepsilon^{-p}|\Omega|.
 \end{align}
Inserting \eqref{ge7}-\eqref{ge9} in \eqref{ge6} yields
\begin{align}
\label{ge10}
  \nonumber \dfrac{d}{dt}\int_\Omega c_{1\varepsilon}^{p} + \frac{p(p - 1)a_1}{2} \int_\Omega c_{1 \varepsilon}^{p - 2}|\nabla c_{1 \varepsilon}|^2 & \leq \dfrac{p(p - 1)}{a_1} b_h^2\int_\Omega c_{1 \varepsilon}^{p} |\nabla h_{\varepsilon}|^2 + \dfrac{p(p - 1)}{a_1} b_\tau^2\int_\Omega c_{1 \varepsilon}^{p} |\nabla \tau_{\varepsilon}|^2  \\[5pt]
  & + p(\beta + 1) \int_\Omega  c_{1 \varepsilon}^{p} + C(p)  M_{\alpha_2}^p \varepsilon^{-p}|\Omega|
\end{align}
for all $t \in (0, T_{\max, \varepsilon})$. Thanks to \eqref{ge3} we obtain
\begin{align}
\label{ge11}
  \nonumber \dfrac{d}{dt}\int_\Omega c_{1\varepsilon}^{p} & + \frac{p(p - 1)a_1}{2} \int_\Omega c_{1 \varepsilon}^{p - 2}|\nabla c_{1 \varepsilon}|^2 \leq \bigg\{\frac{p(p - 1)}{a_1} b_h^2 \|\nabla h_\varepsilon(\cdot, t)\|^2_{L^\infty(\Omega)} \\ & + \frac{p(p - 1)}{a_1} b_\tau^2 \|\nabla \tau_\varepsilon(\cdot, t)\|^2_{L^\infty(\Omega)} + p(\beta + 1) \bigg\}\int_\Omega  c_{1 \varepsilon}^{p} + C(p)  M_{\alpha_2}^p \varepsilon^{-p}|\Omega|
\end{align}
for all $t \in (0, \widehat{T}_{\varepsilon})$. Utilizing Gronwall's inequality in view of \eqref{ge3} we conclude that
\begin{align}
  \label{c1pbound}
  \nonumber \int_\Omega c^p_{1\varepsilon}(\cdot, t) \leq & \exp\bigg\{\frac{p(p - 1)b_h^2}{a_1} \int_0^t \|\nabla h_\varepsilon(\cdot, s)\|_{L^\infty(\Omega)}^2 ds + \frac{p(p - 1)b_\tau^2}{a_1} \int_0^t \|\nabla \tau_\varepsilon(\cdot, s)\|_{L^\infty(\Omega)}^2 ds\\ & + p(\beta + 1)t\bigg\} \left(\int_\Omega c_{10\varepsilon} + C(p)  M_{\alpha_2}^p \varepsilon^{-p}|\Omega| t\right) \leq C_5(\varepsilon, T, p) \quad \text{for all}~ t \in (0, \widehat{T}_\varepsilon).
\end{align}
This allows us to choose a $C_6 = C_6(\varepsilon, T) > 0$ fulfilling
\begin{equation}
\label{c2bound11}
  \|c_{1\varepsilon}(\cdot, t)\|_{L^{ p}(\Omega)} \leq C_6 \quad \text{for all}~ t \in (0, \widehat{T}_{\varepsilon}).
\end{equation}
For all $t \in (0, \widehat{T}_{\varepsilon})$ we can have
\begin{equation}
  \label{c2bound12}
  \bar{c}_{1\varepsilon}(t) := \frac{1}{|\Omega|} \int_\Omega c_{1\varepsilon}(\cdot, t) \leq |\Omega|^{-\frac{1}{p}} \|c_{1\varepsilon}(\cdot, t)\|_{L^{p}(\Omega)} \leq |\Omega|^{-\frac{1}{p}} C_6,
\end{equation}
and
\begin{equation}
  \label{c2bound13}
  \|c_{1\varepsilon}(\cdot, t) - \bar{c}_{1\varepsilon}(t)\|_{L^{p}(\Omega)} \leq 2 C_6.
\end{equation}
Applying the variation of constants formula to the second equation of \eqref{model1} we estimate
\begin{align}
  \label{c2_bound1}
    \nonumber \|c_{2\varepsilon}(\cdot, t)\|_{L^\infty(\Omega)} & = \sup_{\Omega} \big(e^{a_2 t \Delta }c_{20\varepsilon} + \int_0^t e^{(t - s) a_2 \Delta} \alpha_1(\tau_\varepsilon)(\cdot, {s}) c_{1\varepsilon}(\cdot, {s}) ds \\[5pt]
    \nonumber & - \int_0^t e^{(t - s) a_2 \Delta} \alpha_2(\tau_\varepsilon)(\cdot, {s}) F_\varepsilon(c_{2\varepsilon}(\cdot, s)) ds \bigg)   \\[5pt]
    \nonumber & \leq \sup_{\Omega} \left(e^{a_2 t \Delta }c_{20\varepsilon}\right) + \sup_{\Omega} \left( \int_0^t e^{(t - s) a_2 \Delta} \alpha_1(\tau_\varepsilon)(\cdot, s) c_{1\varepsilon}(\cdot, s) ds\right)  \\[5pt]
    \nonumber & \leq\|e^{a_2 t \Delta }c_{20\varepsilon}\|_{L^\infty(\Omega)} + M_{\alpha_1} \int_0^t \|e^{(t - s) a_2 \Delta} c_{1\varepsilon}(\cdot, s)\|_{L^\infty(\Omega)} ds  \\[5pt]
    \nonumber & \leq \|e^{a_2 t \Delta }c_{20\varepsilon}\|_{L^\infty(\Omega)} + M_{\alpha_1}\int_0^t \|  e^{(t - s) a_2 \Delta} (c_{1\varepsilon}(\cdot, s) - \bar{c}_{1\varepsilon}(s))\|_{L^\infty(\Omega)} ds \\[5pt]
    & + M_{\alpha_1}  \int_0^t  \|  e^{(t - s) a_2 \Delta}\bar{c}_{1\varepsilon}(s)\|_{L^\infty(\Omega)} ds \quad  ~\text{for all}~ t \in (0, \widehat{T}_{\varepsilon}).
  \end{align}
 The parabolic maximum principle entails that
 \begin{equation}
   \label{c2bound14}
   \|e^{a_2 t \Delta }c_{20\varepsilon}\|_{L^\infty(\Omega)} \leq C_7
 \end{equation}
and
\begin{equation}
  \label{c2bound15}
  M_{\alpha_1}  \int_0^t  \|  e^{(t - s) a_2 \Delta}\bar{c}_{1\varepsilon}(s) \|_{L^\infty(\Omega)} ds \leq M_{\alpha_1}  \int_0^t \bar{c}_{1\varepsilon}(s)  ds \leq M_{\alpha_1} |\Omega|^{-\frac{1}{n}}C_6 \widehat{T}_\varepsilon
\end{equation}
for all $t \in (0, \widehat{T}_{\varepsilon})$. Applying the properties of the Neumann heat semigroup estimates \cite[Lemma 1.3]{winkler2010aggregation}, we can find a positive constant $C_8 > 0$ such that
\begin{align}
  \nonumber & M_{\alpha_1}\int_0^t \|  e^{(t - s) a_2 \Delta} (c_{1\varepsilon}(\cdot, s) - \bar{c}_{1\varepsilon}(s)\|_{L^\infty(\Omega)} ds  \\[5pt]
  \nonumber & \leq M_{\alpha_1} C_8 \int_0^t (1 + (a_2(t - s))^{-  \frac{n}{2p}})e^{-a_2\lambda_1(t - s)}\|(c_{1\varepsilon}(\cdot, s) - \bar{c}_{1\varepsilon}(s)\|_{L^p(\Omega)} ds  \\[5pt]
  \nonumber & \leq 2 C_6 C_8M _{\alpha_1}\int_0^t (1 + (a_2(t-s)^{- \frac{n}{2p}})e^{-a_2\lambda_1 (t-s)}ds \  \\[5pt]
  & \leq C_9 \quad \text{for all}~ t \in (0, \widehat{T}_{\varepsilon}),
\end{align}
with $C_9 = C_9(\varepsilon, T) > 0$, where in the last inequality we have used the fact that $\int_0^\infty  (1 + \zeta^{- \frac{n}{2p}})e^{-\lambda_1 \zeta }d\zeta   < \infty$ for $p>n/2$.  Here $\lambda_1 > 0$ is the first positive eigenvalue of the Laplacian operator coupled with Neumann boundary conditions in the domain $\Omega$. From above we can see that
\begin{equation}
  \label{c2_bound}
  \|c_{2\varepsilon}(\cdot, t)\|_{L^\infty(\Omega)} \leq C_9(\varepsilon, T).
\end{equation}
\noindent
Because of \eqref{h_bound}, \eqref{t_bound},  \eqref{c2bound11}, \eqref{c2_bound}, and the fact that $\frac{c_{2\varepsilon}}{1 + c_{2\varepsilon}} \leq 1$ for $c_{2\varepsilon} \geq 0$, we can claim that $f_\varepsilon := - \gamma_1 h_\varepsilon c_{1\varepsilon}  - \gamma_2 h_\varepsilon c_{2\varepsilon} - \mu h_\varepsilon +  \tfrac{c_{2\varepsilon}}{1 + c_{2\varepsilon}} $ and $g_\varepsilon :=  - \delta c_{1\varepsilon} \tau_\varepsilon - \sigma \tau_\varepsilon + \frac{c_{2\varepsilon}}{1 + c_{2\varepsilon}} $ belong to $L^\infty((0, \widehat{T}_\varepsilon), L^p(\Omega))$ for any $p >  n/2$. Using this in conjunction with Neumann heat semigroups \cite[Lemma 1.3]{winkler2010aggregation} and the variation of constants formula, we obtain
\begin{align}
\label{nabh_bound}
  \nonumber & \|\nabla h_\varepsilon(\cdot, t)\|_{L^\infty(\Omega)} \leq \|\nabla e^{\varepsilon t \Delta}h_\varepsilon(\cdot, 0)\|_{L^\infty(\Omega)} + \int_0^t \| \nabla e^{\varepsilon (t - s)\Delta }f_\varepsilon  \|_{L^\infty(\Omega)}ds\\
    & \leq C_{10}(1 + (\varepsilon t)^{-\frac{1}{2}})e^{-\lambda_1 \varepsilon t} \|h_{0\varepsilon}\|_{L^\infty(\Omega)}
   + \int_0^t C_{10}(1 + (\varepsilon (t - s))^{-\frac{1}{2} - \frac{n}{2}\frac{1}{p}})e^{- \lambda_1 \varepsilon(t - s)}\|f_\varepsilon \|_{L^p(\Omega)}ds,
\end{align}
and,
\begin{align}
\label{nabtau_bound}
  \nonumber & \|\nabla \tau_\varepsilon(\cdot, t)\|_{L^\infty(\Omega)} \leq \|\nabla e^{\varepsilon t \Delta}\tau_\varepsilon(\cdot, 0)\|_{L^\infty(\Omega)} + \int_0^t \| \nabla e^{\varepsilon (t - s)\Delta }g_\varepsilon  \|_{L^\infty(\Omega)}ds\\
    & \leq C_{11}(1 + (\varepsilon t)^{-\frac{1}{2}})e^{-\lambda_1 \varepsilon t} \|\tau_{0\varepsilon}\|_{L^\infty(\Omega)}
   + \int_0^t C_{11}(1 + (\varepsilon (t - s))^{-\frac{1}{2} - \frac{n}{2}\frac{1}{p}})e^{- \lambda_1 \varepsilon (t - s)}\|g_\varepsilon \|_{L^p(\Omega)}ds,
\end{align}
for all $t \in (0, \widehat{T}_{\varepsilon})$ and all $p \in (n, \infty)$.  Here we need to require $p>n$ instead of $p>n/2$, in order to ensure convergence of the integrals in \eqref{nabh_bound} and \eqref{nabtau_bound}, thereby also using the above mentioned boundedness of $\|f_\varepsilon(\cdot, s)\|_{L^p(\Omega)}$ and $\|g_\varepsilon(\cdot, s)\|_{L^p(\Omega)}$  for $s \in (0, \widehat{T}_{\varepsilon})$.\\[-2ex]

\noindent
These estimates enable us to establish global existence in time for \eqref{model1}. To this end we define the following two functions:
\begin{align*}
  & A(x, t, c_{1\varepsilon}, q) := a_1 q - b_h c_{1\varepsilon} \nabla h_\varepsilon(x, t) - b_\tau c_{1\varepsilon} \nabla \tau_\varepsilon(x, t),          \\[5pt]
  & B(x, t, c_{1\varepsilon}, c_{2\varepsilon}) := - \alpha_1(\tau_\varepsilon)  c_{1\varepsilon}   + \alpha_2(\tau_\varepsilon)  F_\varepsilon (c_{2\varepsilon})
    + \beta c_{1\varepsilon}(1 - c_{1\varepsilon} - c_{2\varepsilon} - \tau_\varepsilon) - \varepsilon c_{1\varepsilon}^\theta,
\end{align*}
in $\Omega \times (0, T_{\max, \varepsilon}) \times \mathbb{R} \times \mathbb{R}^n$ and $\Omega \times (0, T_{\max, \varepsilon}) \times \R\times \R$, respectively. Then \eqref{model1} writes
\begin{align*}
  & \partial_t c_{1\varepsilon} - \nabla \cdot (A(x, t, c_{1\varepsilon}, \nabla c_{1\varepsilon})) = B(x, t, c_{1\varepsilon}, c_{2\varepsilon})
\end{align*}
holding for all $(x, t) \in \Omega \times (0, T_{\max, \varepsilon}) $. We also get
\begin{equation*}
  \begin{cases}
    & A(x, t, c_{1\varepsilon}, \nabla c_{1\varepsilon})\nabla c_{1\varepsilon}  \geq \frac{a_1}{2} |\nabla c_{1\varepsilon}|^2 - \psi_1(x, t),  \\[5pt]
    & |A(x, t, c_{1\varepsilon}, \nabla c_{1\varepsilon})| \leq a_1 |\nabla c_{1\varepsilon}| + \psi_2(x, t),   \\[5pt]
    & |B(x, t, c_{1\varepsilon}, c_{2\varepsilon})| \leq \psi_3(x, t),
  \end{cases}
\end{equation*}
where $\psi_1: = \frac{b_h^2}{a_1} c_{1\varepsilon}^2 |\nabla h_\varepsilon|^2 +  \frac{b_\tau^2}{a_1} c_{1\varepsilon}^2 |\nabla \tau_\varepsilon|^2,\ \psi_2 :=  b_h c_{1\varepsilon} |\nabla h_\varepsilon| + b_\tau c_{1\varepsilon} |\nabla \tau_\varepsilon|,\  \psi_3 := \alpha_1(\tau_\varepsilon) c_{1\varepsilon}  + \alpha_2(\tau_\varepsilon) F_\varepsilon(c_{2\varepsilon}) + \beta c_{1\varepsilon} + \beta c_{1\varepsilon}^2 + \beta c_{1\varepsilon} c_{2\varepsilon} + \beta c_{1\varepsilon} \tau_\varepsilon +\varepsilon c_{1\varepsilon}^\theta$ are nonnegative functions on $\Omega \times (0, T_{\max, \varepsilon})$.\\[-2ex]

\noindent
If we set $T := T_{\max, \varepsilon}$ then by \eqref{c1pbound}, \eqref{c2_bound}, \eqref{nabh_bound},  and \eqref{nabtau_bound} all the $\psi_i's, i = \{1, 2, 3\}$ belong to $L^p(\Omega \times (0, T_{\max, \varepsilon}))$ for every $p > n$. Applying the parabolic H\"{o}lder estimates \cite[Theorem 1.3 and Remark 1.4]{ladyzhenskaia1968linear}, we can claim that there exists a $\vartheta \in (0, 1)$ fulfilling
\begin{equation*}
  \|c_{1\varepsilon}\|_{C^{\vartheta, \frac{\vartheta}{2}}(\Omega \times [0, T_{\max, \varepsilon}])} \leq C_{12}
\end{equation*}
with $C_{12} > 0$. Applying the standard parabolic Schauder estimates \cite{ladyzhenskaia1968linear} to $c_{2\varepsilon}, h_\varepsilon$ and $\tau_\varepsilon$ equations in \eqref{model1}, we can find a $C_{13} > 0$ such that
\begin{equation}
\begin{aligned}
\label{contra1}
  \|c_{2\varepsilon}\|_{C^{2 + \vartheta, 1 + \frac{\vartheta}{2}}(\Omega \times [0, T_{\max, \varepsilon}])} & \leq C_{13}, ~ \|h_{\varepsilon}\|_{C^{2 + \vartheta, 1 + \frac{\vartheta}{2}}(\Omega \times [0, T_{\max, \varepsilon}])} \leq C_{13} ~ \text{and}~\\[5pt]  & \|\tau_{\varepsilon}\|_{C^{2 + \vartheta , 1 + \frac{\vartheta}{2}}(\Omega \times [0, T_{\max, \varepsilon}])} \leq C_{13}.
\end{aligned}
\end{equation}
Again, the standard parabolic Schauder estimates enable us to find a $C_{14} > 0$ such that
\begin{equation}
\label{contra2}
  \|c_{1\varepsilon}\|_{C^{2 + \vartheta, 1 + \frac{\vartheta}{2}}(\Omega \times [0, T_{\max, \varepsilon}])} \leq C_{14}.
\end{equation}
Taken together \eqref{contra1} and \eqref{contra2} contradict the extensibility criterion \eqref{loe1}.

\end{proof}

%%%%%%%%%%%%%%%%%%%%%%%%%%%%%%%%%%%%%%%%%%%%%%%%%%%%%%%%%%%%%%%%%%%%%%%%%%%%%%%%%%%%%%%%%%%%%%%%%%%%%%%%%%%%%%%%%%%%%%%%%%%%%%%%%%%%%%%%%%%%%%%%%%%%%%%%%%%%%%%%%%%%%%%%%%%%%%%%%%%%%%%%%%%%%%%%%%%%%%%%%%%%%%%
%%%%%%%%%%%%%%%%%%%%%%%%%%%%%%%%%%%%%%%%%%%%%%%%%%%%%%%%%%%%%%%%%%%%%%%%%%%%%%%%%%%%%%%%%%%%%%%%%%%%%%%%%%%%%%%%%%%%%%%%%%%%%%%%%%%%%%%%%%%%%%%%%%%%%%%%%%%%%%%%%%%%%%%%%%%%%%%%%%%%%%%%%%%%%%%%%%%%%%%%%%%%%%

\section{An entropy-type functional}

The aim of this section is to derive some estimates which stem from an entropy-type functional. This is the main step towards the existence of a global weak solution to \eqref{model1}. We will initially establish certain inequalities that will prove to be valuable in the sequel.

\begin{Lemma}
\label{elemmac1}
   For all $t > 0$ and for every $\varepsilon \in (0, 1)$, there exist some constants $C,C_1,C_2 > 0$ such that
  \begin{align}
    \label{ey1}
    & \nonumber \dfrac{d}{dt} \int_\Omega c_{1 \varepsilon} \ln c_{1 \varepsilon} + a_1\int_\Omega \dfrac{|\nabla c_{1 \varepsilon}|^2}{c_{1\varepsilon}} + \dfrac{\varepsilon}{2} \int_\Omega c_{1 \varepsilon}^\theta \ln(2 + c_{1 \varepsilon}) + \dfrac{\beta}{2} \int_\Omega c_{1 \varepsilon}^2 \ln(2 + c_{1 \varepsilon}) \\[5pt]
    & \leq b_h \int_\Omega \nabla c_{1 \varepsilon} \cdot \nabla h_\varepsilon + b_\tau \int_\Omega \nabla c_{1 \varepsilon} \cdot \nabla \tau_\varepsilon + M_{\alpha_2} \int_{\{c_{1\varepsilon} > 1\}} c_{2\varepsilon} \ln c_{1\varepsilon} + C_1\int_\Omega c_{1\varepsilon} + C_2 \int_\Omega c_{2\varepsilon} + C.
    \end{align}
\end{Lemma}
\begin{proof}
  The positivity of $c_{1\varepsilon}$ allows us to test the first equation of \eqref{model1} with $\ln c_{1\varepsilon}$:
     \begin{align}
    \label{ey2}
    \nonumber & \dfrac{d}{dt} \int_\Omega c_{1\varepsilon} \ln c_{1\varepsilon} = \int_\Omega \partial_t c_{1\varepsilon} \left(\ln c_{1\varepsilon} + 1 \right) \\
    \nonumber & = - a_1\int_\Omega \dfrac{|\nabla c_{1 \varepsilon}|^2}{c_{1\varepsilon}} + b_h \int_\Omega \nabla  c_{1\varepsilon} \cdot \nabla h_\varepsilon +  b_\tau \int_\Omega \nabla  c_{1\varepsilon} \cdot \nabla \tau_\varepsilon - \int_\Omega \alpha_1(\tau_\varepsilon) c_{1\varepsilon} \ln c_{1\varepsilon} \\
    &\nonumber + \int_\Omega \alpha_2(\tau_\varepsilon) F_\varepsilon(c_{2\varepsilon}) \ln c_{1\varepsilon} + \int_\Omega \beta c_{1\varepsilon}(1 - c_{1\varepsilon} -  c_{2\varepsilon} -  \tau_{\varepsilon})\ln c_{1\varepsilon} - \int_\Omega \varepsilon c_{1\varepsilon}^\theta \ln c_{1\varepsilon} - \int_\Omega \alpha_1(\tau_\varepsilon) c_{1\varepsilon} \\
    &  + \int_\Omega \alpha_2(\tau_\varepsilon) F_\varepsilon(c_{2\varepsilon}) + \int_\Omega \beta c_{1\varepsilon} (1 - c_{1\varepsilon} - c_{2\varepsilon} - \tau_\varepsilon) - \varepsilon\int_\Omega c_{1\varepsilon}^\theta,
  \end{align}
  for all $t > 0$. The first three terms on the right-hand side of (\ref{ey2}) have been obtained using integration by parts.\\[-2ex]
  
  \noindent
  Using \cite[Lemma 4.2]{stinner2014global} we can have
  \begin{align}
    \label{ey3}
    & \int_\Omega \beta(c_{1\varepsilon} \ln c_{1\varepsilon} - c_{1\varepsilon}^2 \ln c_{1\varepsilon}) \leq - \dfrac{\beta}{2} \int_\Omega c_{1\varepsilon}^2 \ln(2 + c_{1\varepsilon}) + C_3 \quad ~\text{and}~\\
    \label{ey4}
    & - \int_\Omega \varepsilon c_{1\varepsilon}^\theta \ln c_{1\varepsilon} \leq - \dfrac{\varepsilon}{2} \int_\Omega c_{1\varepsilon}^\theta \ln(2 + c_{1\varepsilon}) + C_4.
  \end{align}
  Using the nonnegativity of solution components, \eqref{alpha_def}, \eqref{t_bound} and the inequality $s \ln s \geq - \frac{1}{e}$ for all $s > 0$, we can  have
  \begin{align}
    \label{ey5} & -\int_\Omega \alpha_1(\tau_\varepsilon) c_{1\varepsilon} \ln c_{1\varepsilon} = \int_\Omega \alpha_1(\tau_\varepsilon) (- c_{1\varepsilon} \ln c_{1\varepsilon})  \leq  \tfrac{1}{e}M_{\alpha_1} |\Omega|,                                 \\[5pt]
    \label{ney9}& \int_\Omega \alpha_2(\tau_\varepsilon) F_\varepsilon(c_{2\varepsilon}) \ln c_{1\varepsilon} = \int_{\Omega} \alpha_2(\tau_\varepsilon) \tfrac{1}{1 + \varepsilon c_{2\varepsilon}} (c_{2\varepsilon} \ln c_{1\varepsilon}) \leq M_{\alpha_2} \int_{\{c_{1\varepsilon} > 1\}} c_{2\varepsilon} \ln c_{1\varepsilon},       \\[5pt]
    \label{ey6} & - \beta \int_\Omega c_{1\varepsilon} \ln c_{1\varepsilon} \tau_\varepsilon \leq \tfrac{1}{e} \beta M_\tau |\Omega|,                                         \\[5pt]
    \label{ey7} & - \beta \int_\Omega c_{1\varepsilon} \ln c_{1\varepsilon} c_{2\varepsilon} \leq \dfrac{\beta }{e}\int_\Omega c_{2\varepsilon}  \\[5pt]
    \label{ey9}
     \nonumber & -\int_\Omega \alpha_1(\tau_\varepsilon) c_{1\varepsilon} + \int_\Omega \alpha_2(\tau_\varepsilon) F_\varepsilon(c_{2\varepsilon}) + \int_\Omega \beta c_{1\varepsilon} (1 - c_{1\varepsilon} - c_{2\varepsilon} - \tau_\varepsilon) - \varepsilon\int_\Omega c_{1\varepsilon}^\theta \\[5pt] & \leq C_1\int_\Omega c_{1\varepsilon} + C_2 \int_\Omega c_{2\varepsilon},
  \end{align}
with adequate $C_1, C_2 > 0$. Inserting \eqref{ey3}-\eqref{ey9} in \eqref{ey2} we can directly have \eqref{ey1}.

\end{proof}
%%%%%%%%%%%%%%%%%%%%%%%%%%%%%%%%%%%%%%%%%%%%%%%%%%%%%%%%%%%%%%%%%%%%%%%%%%%%%%%%%%%%%%%%%%%%%%%%%%%%%%%%%%%%%%%%%%%%%%%%%%%%%%%%%%%%%%%%%%%%%%

\begin{Lemma}
\label{elemmac2}
 For all $t > 0$ and for every $\varepsilon \in (0, 1)$, there exists a $C > 0$ such that
 \begin{align}
    \label{eyc1}
    \dfrac{d}{dt} \int_\Omega c_{2 \varepsilon} \ln c_{2 \varepsilon} + a_2 \int_\Omega \dfrac{|\nabla c_{2 \varepsilon}|^2}{c_{2\varepsilon}} \leq M_{\alpha_1} \int_{\{c_{2\varepsilon} > 1\}} c_{1\varepsilon} \ln c_{2\varepsilon} + M_{\alpha_1} \int_\Omega c_{1\varepsilon} + M_{\alpha_2} \int_\Omega c_{2\varepsilon} + C.
  \end{align}
  \end{Lemma}

\begin{proof}
  We test the second equation of (\ref{model1}) against $\ln c_{2\varepsilon}$,
   \begin{align}
    \label{eyc2}
    \dfrac{d}{dt} \int_\Omega c_{2 \varepsilon} \ln c_{2 \varepsilon}  & = - a_2 \int_\Omega \dfrac{|\nabla c_{2 \varepsilon}|^2}{c_{2\varepsilon}} + \int_\Omega \alpha_1 (\tau_\varepsilon) c_{1\varepsilon} \ln c_{2\varepsilon} - \int_\Omega \alpha_2(\tau_\varepsilon) F_\varepsilon (c_{2\varepsilon}) \ln c_{2\varepsilon}+ \int_\Omega \alpha_1(\tau_\varepsilon) c_{1\varepsilon} - \int_\Omega \alpha_2(\tau_\varepsilon) F_\varepsilon(c_{2\varepsilon}) ,
   \end{align}
  for all $t > 0$.
  Using the nonnegativity of solution components, \eqref{alpha_def} and the identity $s \ln s \geq - \frac{1}{e}$ for all $s > 0$, we can  have
  \begin{align}
    \label{c2ent1new45} & \int_\Omega \alpha_1(\tau_\varepsilon) c_{1\varepsilon} \ln c_{2\varepsilon} \leq M_{\alpha_1} \int_{\{c_{2\varepsilon} > 1\}} c_{1\varepsilon} \ln c_{2\varepsilon}, \\[5pt]
    \label{c2ent1}& - \int_\Omega \alpha_{2} (\tau_\varepsilon) F_\varepsilon(c_{2 \varepsilon}) \ln c_{2 \varepsilon} = \int_\Omega \alpha_2(\tau_\varepsilon)  \tfrac{1}{1 + \varepsilon c_{2\varepsilon}} (- c_{2\varepsilon} \ln c_{2\varepsilon})  \leq  \tfrac{1}{e}M_{\alpha_2} |\Omega|, \\[5pt]
    \label{c2ent2} & \int_\Omega \alpha_1(\tau_\varepsilon) c_{1\varepsilon} - \int_\Omega \alpha_2(\tau_\varepsilon)  F_\varepsilon(c_{2\varepsilon})  \leq M_{\alpha_1} \int_\Omega c_{1\varepsilon} + M_{\alpha_2} \int_\Omega c_{2\varepsilon}.
  \end{align}
   Inserting \eqref{c2ent1new45}-\eqref{c2ent2} in \eqref{eyc2} we can directly have \eqref{eyc1}.

\end{proof}
%%%%%%%%%%%%%%%%%%%%%%%%%%%%%%%%%%%%%%%%%%%%%%%%%%%%%%%%%%%%%%%%%%%%%%%%%%%%%%%%%%%%%%%%%%%%%%%%%%%%%%%%%%%%%%%%%%%%%%%%%%%%%%%%%%%%%%%%%%%%%%
%%%%%%%%%%%%%%%%%%%%%%%%%%%%%%%%%%%%%%%%%%%%%%%%%%%%%%%%%%%%%%%%%%%%%%%%%%%%%%%%%%%%%%%%%%%%%%%%%%%%%%%%%%%%%%%%%%%%%%%%%%%%%%%%%%%%%%%%%%%%%%

\begin{Lemma}
  \label{elemmah}
   For all $t > 0$ and for every $\varepsilon \in (0, 1)$, there exists a $C > 0$ such that
  \begin{equation}
    \label{eh1}
    \dfrac{1}{2} \dfrac{d}{dt}\int_\Omega \dfrac{|\nabla h_\varepsilon|^2}{h_\varepsilon} + \dfrac{\mu}{2} \int_\Omega \dfrac{|\nabla h_\varepsilon|^2}{h_\varepsilon} + \gamma_1 \int_\Omega \nabla c_{1\varepsilon} \cdot \nabla h_\varepsilon + \dfrac{\gamma_1}{2} \int_\Omega \dfrac{|\nabla h_\varepsilon|^2}{h_\varepsilon} c_{1\varepsilon} \leq (4 \gamma_2 M_h + 1) \int_\Omega  \dfrac{|\nabla c_{2 \varepsilon}|^2}{c_{2\varepsilon}} + C.
  \end{equation}
\end{Lemma}
\begin{proof}

  We can directly compute
  \begin{align}
  \label{eh2}
    & \nonumber  \dfrac{1}{2} \dfrac{d}{dt}\int_\Omega \dfrac{|\nabla h_\varepsilon|^2}{h_\varepsilon}  = \int_\Omega \dfrac{\nabla h_\varepsilon \cdot \nabla h_{\varepsilon t}}{h_\varepsilon} - \dfrac{1}{2} \int_\Omega \dfrac{|\nabla h_\varepsilon|^2}{h_\varepsilon^2} h_{\varepsilon t}                                                                                                                                 \\[7pt]
    & \nonumber = \int_\Omega \dfrac{\nabla h_\varepsilon}{h_\varepsilon} \cdot \nabla (\varepsilon \Delta h_\varepsilon - \gamma_1 h_\varepsilon c_{1\varepsilon}  - \gamma_2 h_\varepsilon c_{2\varepsilon} - \mu h_\varepsilon +  \tfrac{c_{2\varepsilon}}{1 + c_{2\varepsilon}} ) \\[5pt]
    & \nonumber - \dfrac{1}{2} \int_\Omega \dfrac{|\nabla h_\varepsilon|^2}{h_\varepsilon^2} (\varepsilon \Delta h_\varepsilon - \gamma_1 h_\varepsilon c_{1\varepsilon} - \gamma_2 h_\varepsilon c_{2\varepsilon} - \mu h_\varepsilon +  \tfrac{c_{2\varepsilon}}{1 + c_{2\varepsilon}} )\\[7pt]
    \nonumber& = \varepsilon \int_\Omega \tfrac{\nabla h_\varepsilon \cdot \nabla \Delta h_\varepsilon}{h_\varepsilon} - \tfrac{\varepsilon}{2} \int_\Omega \tfrac{|\nabla h_\varepsilon|^2}{h_\varepsilon^2} \Delta h_\varepsilon - \gamma_1 \int_\Omega \nabla c_{1 \varepsilon} \cdot \nabla h_\varepsilon - \tfrac{\gamma_1}{2} \int_\Omega \tfrac{|\nabla h_\varepsilon|^2}{h_\varepsilon} c_{1\varepsilon} - \gamma_2 \int_\Omega  \nabla c_{2 \varepsilon} \cdot \nabla h_\varepsilon \\[7pt]
    & - \frac{\gamma_2}{2} \int_\Omega \frac{|\nabla h_\varepsilon|^2}{h_\varepsilon} c_{2\varepsilon} - \frac{\mu}{2} \int_\Omega \frac{|\nabla h_\varepsilon|^2}{h_\varepsilon} + \int_\Omega \frac{\nabla c_{2\varepsilon} \cdot \nabla h_\varepsilon}{h_\varepsilon (1 + c_{2\varepsilon})^2} - \frac{1}{2} \int_\Omega \frac{c_{2\varepsilon}}{1 + c_{2\varepsilon}} \frac{|\nabla h_\varepsilon|^2}{h_\varepsilon^2} \quad ~\text{for all}~ t > 0.
  \end{align}

  We will now estimate the terms on the right-hand side of (\ref{eh2}). To handle the first two terms we will adopt the approach of \cite[Lemma 3.4]{tao2021global},
  \begin{equation*}
  \begin{aligned}
   & \int_\Omega \dfrac{1}{h_\varepsilon} \nabla \Delta h_\varepsilon \cdot \nabla h_\varepsilon - \dfrac{1}{2}\int_\Omega \dfrac{|\nabla h_\varepsilon|^2}{h_\varepsilon^2} \Delta h_\varepsilon               \\[7pt]
   & = \sum_{i, j = 1}^n \int_\Omega \dfrac{1}{h_\varepsilon} \partial_i h_\varepsilon \partial_{ijj} h_\varepsilon - \dfrac{1}{2} \sum_{i, j = 1}^n \int_\Omega \dfrac{1}{h_\varepsilon^2}\left(\partial_i h_\varepsilon\right)^2 \partial_{jj} h_\varepsilon                                                                                                                                               \\[7pt]
   & = 2 \sum_{i, j = 1}^n \int_\Omega \dfrac{1}{h_\varepsilon^2} \partial_i h_\varepsilon \partial_j h_\varepsilon \partial_{ij} h_\varepsilon - \sum_{i, j = 1}^n \int_\Omega \dfrac{1}{h_\varepsilon} \left(\partial_{ij} h_\varepsilon\right)^2 - \sum_{i, j = 1}^n \int_\Omega \dfrac{1}{h_\varepsilon^3} \left(\partial_i h_\varepsilon\right)^2 \left(\partial_j h_\varepsilon \right)^2                               \\[7pt]
   & + \sum_{i, j = 1}^n \int_{\partial \Omega} \dfrac{1}{h_\varepsilon} \partial_i h_\varepsilon \partial_{ij} h_\varepsilon \nu^j - \dfrac{1}{2} \sum_{i, j = 1}^n \int_{\partial \Omega} \dfrac{1}{h_\varepsilon^2}\left(\partial_i h_\varepsilon\right)^2 \partial_j h_\varepsilon \nu^j, ~~\text{for all} ~~ t > 0,
   \end{aligned}
  \end{equation*}
  where $\nu^j$ is the $j^{th}$ component of $\nu$. A direct computation yields:
  \begin{align}
      \nonumber & 2 \sum_{i, j = 1}^n \int_\Omega \dfrac{1}{h_\varepsilon^2} \partial_i h_\varepsilon \partial_j h_\varepsilon \partial_{ij} h_\varepsilon - \sum_{i, j = 1}^n \int_\Omega \dfrac{1}{h_\varepsilon} \left(\partial_{ij} h_\varepsilon\right)^2 - \sum_{i, j = 1}^n \int_\Omega \dfrac{1}{h_\varepsilon^3} \left(\partial_i h_\varepsilon\right)^2 \left(\partial_j h_\varepsilon \right)^2                               \\[7pt]
      \label{nu1}& = - \sum_{i, j = 1}^n \int_\Omega \dfrac{1}{h_\varepsilon} \left|\dfrac{1}{h_\varepsilon} \partial_i h_\varepsilon \partial_j h_\varepsilon - \partial_{ij} h_\varepsilon \right|^2
      \leq 0 \quad ~\text{for all}~ t > 0.
    \end{align}
  As $\Omega$ is a convex domain and the boundary condition of $h_\varepsilon$ is of the zero-flux type we have (cf. \cite[Lemma 3.2 ]{tao2012boundedness})
    \begin{align}
     \nonumber & \sum_{i, j = 1}^n \int_{\partial \Omega} \dfrac{1}{h_\varepsilon} \partial_i h_\varepsilon \partial_{ij} h_\varepsilon \nu^j = \dfrac{1}{2} \int_{\partial\Omega} \dfrac{1}{h_\varepsilon} \dfrac{\partial |\nabla h_\varepsilon|^2}{\partial \nu} \leq 0 ~~ \text{and},     \\
     \label{nu2}& - \dfrac{1}{2} \sum_{i, j = 1}^n \int_{\partial \Omega} \dfrac{1}{h_\varepsilon^2}\left(\partial_i h_\varepsilon\right)^2 \partial_j h_\varepsilon \nu^j = - \dfrac{1}{2} \int_{\partial\Omega} \dfrac{|\nabla h_\varepsilon|^2}{h_\varepsilon^2} \dfrac{\partial h_\varepsilon}{\partial \nu} = 0 \quad ~\text{for all}~ t > 0.
      \end{align}
    \eqref{nu1} and \eqref{nu2} together result in
    \begin{equation*}
      \varepsilon \int_\Omega \dfrac{1}{h_\varepsilon} \nabla \Delta h_\varepsilon \cdot \nabla h_\varepsilon - \dfrac{\varepsilon}{2}\int_\Omega \dfrac{|\nabla h_\varepsilon|^2}{h_\varepsilon^2} \Delta h_\varepsilon \leq 0  \quad ~\text{for all}~ t > 0.
    \end{equation*}
  Now consider the term $- \gamma_2 \int_\Omega \nabla h_\varepsilon \cdot \nabla c_{2 \varepsilon}$. An application of Young's inequality in conjunction with \eqref{h_bound} yields
  \begin{align}
    \label{neweh1} \nonumber - \gamma_2 \int_\Omega \nabla h_\varepsilon \cdot \nabla c_{2 \varepsilon} & \leq \dfrac{\gamma_2}{2} \int_\Omega \dfrac{|\nabla h_\varepsilon|^2}{h_\varepsilon} c_{2\varepsilon} + 4 \gamma_2 \int_\Omega \dfrac{|\nabla c_{2\varepsilon}|^2}{c_{2 \varepsilon}} h_\varepsilon   \\[7pt]
    & \leq \dfrac{\gamma_2}{2} \int_\Omega \dfrac{|\nabla h_\varepsilon|^2}{h_\varepsilon} c_{2\varepsilon} + 4 \gamma_2 M_h \int_\Omega \dfrac{|\nabla c_{2\varepsilon}|^2}{c_{2 \varepsilon}}
  \end{align}
  for all $t > 0$. To handle the next pair $- \frac{1}{2} \int_\Omega \frac{|\nabla h_\varepsilon|^2}{h_\varepsilon^2} \frac{c_{2\varepsilon}}{1 + c_{2\varepsilon}} + \int_\Omega \frac{\nabla h_\varepsilon \cdot \nabla c_{2\varepsilon}}{h_\varepsilon (1 + c_{2\varepsilon})^2} $ we will proceed as for the estimate derived in (3.13) \cite{surulescu2021does}. It is evident that $\frac{c_{2\varepsilon}}{1 + c_{2\varepsilon}} < 1$ and $\frac{1}{(1 + c_{2\varepsilon})^3} \leq 1$ for $c_{2\varepsilon} \geq 0$, hence
\begin{align}
  \label{eh3}
 - \dfrac{1}{2} \int_\Omega \dfrac{|\nabla h_\varepsilon|^2}{h_\varepsilon^2} \dfrac{c_{2\varepsilon}}{1 + c_{2\varepsilon}} + \int_\Omega \dfrac{\nabla h_\varepsilon \cdot \nabla c_{2\varepsilon}}{h_\varepsilon (1 + c_{2\varepsilon})^2}
%    & \nonumber \leq - \dfrac{1}{2} \int_\Omega \dfrac{|\nabla h_\varepsilon|^2}{h_\varepsilon^2} \dfrac{c_{2\varepsilon}}{1 + c_{2\varepsilon}} + \dfrac{1}{2} \int_\Omega \dfrac{|\nabla h_\varepsilon|^2}{h_\varepsilon^2} \dfrac{c_{2\varepsilon}}{1 + c_{2\varepsilon}} + \dfrac{1}{2} \int_\Omega\dfrac{1}{(1 + c_{2\varepsilon})^4} \cdot \dfrac{1 + c_{2\varepsilon}}{c_{2\varepsilon}} |\nabla c_{2\varepsilon}|^2                                              \\[5pt]
 % & \nonumber =
  & \le \dfrac{1}{2} \int_\Omega\dfrac{1}{(1 + c_{2\varepsilon})^3} \dfrac{|\nabla c_{2\varepsilon}|^2}{c_{2\varepsilon}} \notag \\[5pt]
  & \leq  \int_\Omega \dfrac{|\nabla c_{2\varepsilon}|^2}{c_{2\varepsilon}}\quad \text{for all }t>0.
  \end{align}
%for all $t > 0$.
We get \eqref{eh1} by inserting \eqref{nu1}-\eqref{eh3} into \eqref{eh2}.
\end{proof}
%%%%%%%%%%%%%%%%%%%%%%%%%%%%%%%%%%%%%%%%%%%%%%%%%%%%%%%%%%%%%%%%%%%%%%%%%%%%%%%%%%%%%%%%%%%%%%%%%%%%%%%%%%%%%%%%%%%%%%%%%%%%%%%%%%%%%%%%%%%%%%%%%%%%%%%%%%%%%%%%%%%%%%%%%%%%%%%%%%%%%%%%%%%%%%%%%%%%%%%%%%%%%%

\noindent
Analogously, we can have the following lemma.

%%%%%%%%%%%%%%%%%%%%%%%%%%%%%%%%%%%%%%%%%%%%%%%%%%%%%%%%%%%%%%%%%%%%%%%%%%%%%%%%%%%%%%%%%%%%%%%%%%%%%%%%%%%%%%%%%%%%%%%%%%%%%%%%%%%%%%%%%%%%%%%%%%%%%%%%%%%%%%%%%%%%%%%%%%%%%%%%%%%%%%%%%%%%%%%%%%%%%%%%%%%%%%
\begin{Lemma}
  \label{elemmat}
  For all $t > 0$ and for every $\varepsilon \in (0, 1)$, there exists a $C > 0$ such that
  \begin{equation}
    \label{et1}
    \dfrac{1}{2} \dfrac{d}{dt}\int_\Omega \dfrac{|\nabla \tau_\varepsilon|^2}{\tau_\varepsilon} + \dfrac{\sigma}{2} \int_\Omega \dfrac{|\nabla \tau_\varepsilon|^2}{\tau_\varepsilon} + \delta \int_\Omega \nabla c_{1\varepsilon} \cdot \nabla \tau_\varepsilon + \dfrac{\delta}{2} \int_\Omega \dfrac{|\nabla \tau_\varepsilon|^2}{\tau_\varepsilon} c_{1\varepsilon} \leq  \int_\Omega  \dfrac{|\nabla c_{2 \varepsilon}|^2}{c_{2\varepsilon}} + C.
  \end{equation}
\end{Lemma}

%%%%%%%%%%%%%%%%%%%%%%%%%%%%%%%%%%%%%%%%%%%%%%%%%%%%%%%%%%%%%%%%%%%%%%%%%%%%%%%%%%%%%%%%%%%%%%%%%%%%%%%%%%%%%%%%%%%%%%%%%%%%%%%%%%%%%%%%%%%%%%%%%%%%%%%%%%%%%%%%%%%%%%%%%%%%%%%%%%%%%%%%%%%%%%%%%%%%%%%%%%%%%%
\noindent
Yet other estimates will be needed in order to establish further below some compactness properties of the approximating sequence $(c_{1\varepsilon},c_{2\varepsilon},h_{\varepsilon},\tau _{\varepsilon})_{\varepsilon}$:

\begin{Lemma}
\label{entropy_lemma13}
Let $T > 0$. Then there exists a $C(T) > 0$ such that for all $\varepsilon \in (0, 1)$ the solution to \eqref{model1} satisfies
  \begin{align}
    \label{entropy_main1}
    \int_\Omega  \left(c_{1\varepsilon} \ln c_{1\varepsilon} + \frac{1}{e}\right) + \xi \int_\Omega \left(c_{2\varepsilon} \ln c_{2\varepsilon} + \frac{1}{e}\right) + \dfrac{b_h}{2 \gamma_1} \int_\Omega \dfrac{|\nabla h_\varepsilon|^2}{h_\varepsilon} + \dfrac{b_\tau}{2 \delta} \int_\Omega \dfrac{|\nabla \tau_\varepsilon|^2}{\tau_\varepsilon} \leq C(T) \quad ~\text{for all} ~ t \in (0, T)
     \end{align}
     where $\xi = \tfrac{2}{a_2}((4 \gamma_2 M_h + 1)\tfrac{b_h}{\gamma_1} + \tfrac{b_\tau}{\delta} + 1)$,
     and
    \begin{align}
    \label{entropy_main2}
    \nonumber & \frac{a_1}{4} \int_0^T \int_\Omega \dfrac{|\nabla c_{1 \varepsilon}|^2}{c_{1\varepsilon}} + \frac{1}{2} \int_0^T  \int_\Omega \dfrac{|\nabla c_{2 \varepsilon}|^2}{c_{2\varepsilon}} + \dfrac{\varepsilon}{2} \int_0^T \int_\Omega c_{1 \varepsilon}^\theta \ln(2 + c_{1 \varepsilon}) + \dfrac{\beta}{2} \int_0^T \int_\Omega c_{1 \varepsilon}^2 \ln(2 + c_{1 \varepsilon}) \\[5pt]
  & + \frac{b_h}{2} \int_0^T \int_\Omega \dfrac{|\nabla h_\varepsilon|^2}{h_\varepsilon} c_{1\varepsilon} + \frac{b_\tau}{2} \int_0^T\int_\Omega \dfrac{|\nabla \tau_\varepsilon|^2}{\tau_\varepsilon} c_{1\varepsilon}\leq C(T).
  \end{align}
\end{Lemma}

  \begin{proof}
    Combining \eqref{L1bound} with \eqref{ey1} yields
    \begin{align}
    \label{entropy_new_proof1}
      \nonumber & \dfrac{d}{dt} \int_\Omega \left(c_{1\varepsilon} \ln c_{1\varepsilon} + \frac{1}{e}\right) + \varrho_1 \int_\Omega  \left(c_{1\varepsilon} \ln c_{1\varepsilon} + \frac{1}{e}\right) + a_1\int_\Omega \dfrac{|\nabla c_{1 \varepsilon}|^2}{c_{1\varepsilon}} + \dfrac{\varepsilon}{2} \int_\Omega c_{1 \varepsilon}^\theta \ln(2 + c_{1 \varepsilon}) + \dfrac{\beta}{2} \int_\Omega c_{1 \varepsilon}^2 \ln(2 + c_{1 \varepsilon}) \\[5pt]
      & \leq b_h \int_\Omega \nabla c_{1 \varepsilon} \cdot \nabla h_\varepsilon + b_\tau \int_\Omega \nabla c_{1 \varepsilon} \cdot \nabla \tau_\varepsilon + M_{\alpha_2} \int_{\{c_{1\varepsilon} > 1\}} c_{2\varepsilon} \ln c_{1\varepsilon} + \varrho_1 \int_\Omega c_{1\varepsilon} \ln c_{1\varepsilon} + C_1 C_3 + C_2 C_3 + C + \frac{\varrho_1 |\Omega|}{e}
    \end{align}
    for any $\varrho_1 > 0$ and for all $t \in (0, T)$. To estimate $\varrho_1 \int_\Omega c_{1\varepsilon} \ln c_{1\varepsilon}$ on the right-hand side we will use the fact that $\varrho_1 \int_\Omega c_{1\varepsilon} \ln c_{1\varepsilon} \leq \varrho_1 \int_\Omega c_{1\varepsilon}^{\frac{3}{2}}$, for $c_{1\varepsilon} \geq 0$. Since $n \leq 3$ making use of the Gagliardo-Nirenberg inequality, we employ \eqref{L1bound}  and \cite[(3.19) with $q = \frac{3}{2}$]{tao2021global} to have \begin{align}
    \label{entropy_new_proof2}
        \nonumber \varrho_1 \int_\Omega c_{1\varepsilon}^{\frac{3}{2}} & \leq C_{4} \varrho_1\|\nabla \sqrt{c_{1\varepsilon}}\|^{n/2}_{L^2(\Omega)} \|\sqrt{c_{1\varepsilon}}\|^{3-\frac{n}{2}}_{L^2(\Omega)} + C_4 \varrho_1 \|\sqrt{c_{1\varepsilon}}\|^3_{L^2(\Omega)} \\[5pt]
        & \leq \frac{C_{5} \varrho_1}{4} \int_\Omega  \dfrac{|\nabla c_{1 \varepsilon}|^2}{c_{1\varepsilon}} + C_6 \quad \text{for all} ~ t \in (0, T ),
    \end{align}
    where $C_5 > 0$ and $C_6 > 0$. Setting $\varrho_1 = \frac{2 a_1}{C_5}$ in \eqref{entropy_new_proof2} results in the following estimate
    \begin{align}
    \label{entropy_new_proof3}
        \varrho_1 \int_\Omega c_{1\varepsilon}^{\frac{3}{2}} & \leq  \frac{a_1}{2} \int_\Omega  \dfrac{|\nabla c_{1 \varepsilon}|^2}{c_{1\varepsilon}} + C_6 \quad \text{for all} ~ t \in (0, T ).
    \end{align}
  Using the above estimates in \eqref{entropy_new_proof1} results in
   \begin{align}
    \label{entropy_new_proof4}
      \nonumber & \dfrac{d}{dt} \int_\Omega \left(c_{1\varepsilon} \ln c_{1\varepsilon} + \frac{1}{e}\right) + \varrho_1 \int_\Omega  \left(c_{1\varepsilon} \ln c_{1\varepsilon} + \frac{1}{e}\right) + \frac{a_1}{2}\int_\Omega \dfrac{|\nabla c_{1 \varepsilon}|^2}{c_{1\varepsilon}} + \dfrac{\varepsilon}{2} \int_\Omega c_{1 \varepsilon}^\theta \ln(2 + c_{1 \varepsilon})  \\[5pt]
      & + \dfrac{\beta}{2} \int_\Omega c_{1 \varepsilon}^2 \ln(2 + c_{1 \varepsilon})\leq b_h \int_\Omega \nabla c_{1 \varepsilon} \cdot \nabla h_\varepsilon + b_\tau \int_\Omega \nabla c_{1 \varepsilon} \cdot \nabla \tau_\varepsilon + M_{\alpha_2} \int_{\{c_{1\varepsilon} > 1\}} c_{2\varepsilon} \ln c_{1\varepsilon} + C_7
    \end{align}
    for all $t \in (0, T)$. Following the same procedure we can derive the following estimate from \eqref{eyc1}
    \begin{align}
    \label{entropy_new_proof5}
      \dfrac{d}{dt} \int_\Omega \left(c_{2 \varepsilon} \ln c_{2 \varepsilon} + \frac{1}{e}\right) + \varrho_2  \int_\Omega \left(c_{2 \varepsilon} \ln c_{2 \varepsilon} + \frac{1}{e}\right)  + \frac{a_2}{2} \int_\Omega \dfrac{|\nabla c_{2 \varepsilon}|^2}{c_{2\varepsilon}} \leq M_{\alpha_1} \int_{\{c_{2\varepsilon} > 1\}} c_{1\varepsilon} \ln c_{2\varepsilon}  + C_8
    \end{align}
    for all $t \in (0, T)$. Using \eqref{eh1}, \eqref{et1}, \eqref{entropy_new_proof4} and \eqref{entropy_new_proof5} we can have
  \begin{align}
   \label{entro1}
    \nonumber & \dfrac{d}{dt} \left\{ \int_\Omega \left(c_{1 \varepsilon} \ln c_{1 \varepsilon} + \frac{1}{e}\right) + \tfrac{2}{a_2}((4 \gamma_2 M_h + 1)\tfrac{b_h}{\gamma_1} + \tfrac{b_\tau}{\delta} + 1)\int_\Omega \left(c_{2 \varepsilon} \ln c_{2 \varepsilon} + \frac{1}{e}\right) + \tfrac{b_h}{2 \gamma_1} \int_\Omega \dfrac{|\nabla h_\varepsilon|^2}{h_\varepsilon} + \tfrac{b_\tau}{2 \delta} \int_\Omega \dfrac{|\nabla \tau_\varepsilon|^2}{\tau_\varepsilon}\right\} \\[5pt]
    \nonumber &+ \varrho_1 \int_\Omega \left(c_{1 \varepsilon} \ln c_{1 \varepsilon} + \frac{1}{e}\right) + \tfrac{2}{a_2}\varrho_2((4 \gamma_2 M_h + 1)\tfrac{b_h}{\gamma_1} + \tfrac{b_\tau}{\delta} + 1)\int_\Omega \left(c_{2 \varepsilon} \ln c_{2 \varepsilon} + \frac{1}{e}\right) + \tfrac{b_h \mu }{2 \gamma_1} \int_\Omega \dfrac{|\nabla h_\varepsilon|^2}{h_\varepsilon} + \tfrac{b_\tau \sigma}{2 \delta} \int_\Omega \dfrac{|\nabla \tau_\varepsilon|^2}{\tau_\varepsilon} \\[5pt]
    \nonumber & + \frac{a_1}{2}\int_\Omega \dfrac{|\nabla c_{1 \varepsilon}|^2}{c_{1\varepsilon}} +  ((4 \gamma_2 M_h + 1)\tfrac{b_h}{\gamma_1} + \tfrac{b_\tau}{\delta} + 1)\int_\Omega \dfrac{|\nabla c_{2 \varepsilon}|^2}{c_{2\varepsilon}} + \dfrac{\varepsilon}{2} \int_\Omega c_{1 \varepsilon}^\theta \ln(2 + c_{1 \varepsilon}) \\[5pt]
    \nonumber & + \dfrac{\beta}{2} \int_\Omega c_{1 \varepsilon}^2 \ln(2 + c_{1 \varepsilon}) + \frac{b_h}{2}\int_\Omega \dfrac{|\nabla h_\varepsilon|^2}{h_\varepsilon} c_{1\varepsilon} + \frac{b_\tau}{2} \int_\Omega \dfrac{|\nabla \tau_\varepsilon|^2}{\tau_\varepsilon} c_{1\varepsilon} \leq M_{\alpha_2} \int_{\{c_{1\varepsilon} > 1\}} c_{2\varepsilon} \ln c_{1\varepsilon} \\[5pt]
    &+   \tfrac{2}{a_2}M_{\alpha_1}((4 \gamma_2 M_h + 1)\tfrac{b_h}{\gamma_1} + \tfrac{b_\tau}{\delta} + 1)\int_{\{c_{2 \varepsilon} > 1\}} c_{1\varepsilon} \ln c_{2\varepsilon} + (\tfrac{b_h}{\gamma_1} (4 \gamma_2 M_h + 1) + \tfrac{b_\tau}{\delta})\int_\Omega  \dfrac{|\nabla c_{2 \varepsilon}|^2}{c_{2\varepsilon}} + C_9
  \end{align}
  for all $t \in (0, T)$. As the taxis terms have been counteracted using favorable cancelations, we will now proceed to show that the ill-signed integrals on the right-hand side of \eqref{entro1} can be controlled by the dissipative rates of $\int_\Omega \frac{|\nabla c_{1 \varepsilon}|^2}{c_{1\varepsilon}}$ and $\int_\Omega \frac{|\nabla c_{2 \varepsilon}|^2}{c_{2\varepsilon}}$.\\[-2ex]

  \noindent
  Consider the term $ M_{\alpha_{2}}\int_\Omega c_{2 \varepsilon} \ln c_{1 \varepsilon}$. Two applications of Young's inequality result in two constants $\epsilon_1 > 0$ and $\epsilon_2 > 0$ such that
 \begin{align}
 \label{ent4}
   \nonumber  M_{\alpha_{2}}\int_{\{c_{1\varepsilon} > 1\}} c_{2 \varepsilon} \ln c_{1 \varepsilon} & \leq \epsilon_1 \int_{\{c_{1\varepsilon} > 1\}} c_{2 \varepsilon}^{\frac{5}{3}} + C(\epsilon_1)\left( M_{\alpha_{2}}\right)^{\frac{5}{2}} \int_{\{c_{1 \varepsilon} > 1\}}(\ln c_{1 \varepsilon})^{\frac{5}{2}} \\[7pt]
   & \leq \epsilon_1\int_\Omega c_{2 \varepsilon}^{\frac{5}{3}} + \epsilon_2\left(  M_{\alpha_{2}}\right)^{4} \int_{\{c_{1 \varepsilon} > 1\}} (\ln c_{1 \varepsilon})^{4} + C(\epsilon_2) |\Omega| \quad  ~~\text{for all} ~ t \in (0, T).
 \end{align}
Making use of the Gagliardo-Nirenberg inequality, we employ \eqref{L1bound}  and \cite[(3.19) with $q = \frac{5}{3}$]{tao2021global}  to have
  \begin{align}
   \label{ent5}
   \nonumber \int_\Omega c_{2 \varepsilon}^{\frac{5}{3}} & \leq C_{10}\left(\|\nabla \sqrt{c_{2\varepsilon}}\|_{L^2(\Omega)}^{\frac{2n}{3}} \|\sqrt{c_{2\varepsilon}}\|_{L^2(\Omega)}^{\frac{10}{3} - \frac{5n}{3} + n} + \|\sqrt{c_{2\varepsilon}}\|_{L^2(\Omega)}^{\frac{10}{3}}\right)    \\[7pt]
   & \leq C_{11} \int_\Omega \dfrac{|\nabla c_{2 \varepsilon}|^2}{c_{2\varepsilon}} + C_{12} \quad  ~~\text{for all} ~ t \in (0, T).
 \end{align}
  Using \eqref{L1bound} and the inequality $\ln s \leq \left(\frac{4}{e}\right) s^{\frac{1}{4}}$,  we can have
 \begin{align}
   \label{ent51}
   \epsilon_2\left(  M_{\alpha_{2}}\right)^{4}\int_{\{c_{1 \varepsilon} > 1\}} (\ln c_{1 \varepsilon})^{4} \leq \epsilon_2\left(  M_{\alpha_{2}}\right)^{4} \left(\dfrac{4}{e}\right)^4 \int_\Omega c_{1 \varepsilon}
   \leq C_{13} \quad \text{for all} ~ t \in (0, T).
 \end{align}
  Substituting \eqref{ent5}-\eqref{ent51} into \eqref{ent4} with $\epsilon_1 := \frac{1}{2 C_{11}} > 0$ yields
\begin{align}
  \label{pent1}
  M_{\alpha_{2}}\int_{\{c_{1 \varepsilon} > 1\}} c_{2 \varepsilon} \ln c_{1 \varepsilon} \leq \dfrac{1}{2}  \int_\Omega \dfrac{|\nabla c_{2 \varepsilon}|^2}{c_{2\varepsilon}} + C_{14} \quad \text{for all} ~ t \in (0, T).
\end{align}
Now, consider the term $\tfrac{2}{a_2}[(4 \gamma_2 M_h + 1)\frac{b_h}{\gamma_1} + \frac{b_\tau}{\delta} + 1]M_{\alpha_1} \int_{\{c_{2\varepsilon} > 1\}} c_{1\varepsilon} \ln c_{2\varepsilon}$. A treatment similar to the one above can yield a $C_{15} > 0$ such that %{\cmg Again, I do not believe you can get $a_1/4$ in front of the integral on the RHS.}
\begin{align}
  \label{pent6}
  \tfrac{2}{a_2}M_{\alpha_1} ((4 \gamma_2 M_h + 1)\tfrac{b_h}{\gamma_1} + \tfrac{b_\tau}{\delta} + 1)\int_{\{c_{2\varepsilon} > 1\}} c_{1\varepsilon} \ln c_{2\varepsilon} \leq  \dfrac{a_1}{4} \int_\Omega \dfrac{|\nabla c_{1 \varepsilon}|^2}{c_{1\varepsilon}} +C_{15} \quad \text{for all} ~ t \in (0, T).
\end{align}
Let
\begin{equation*}
\mathcal{F}_\varepsilon(t) := \int_\Omega \left(c_{1 \varepsilon} \ln c_{1 \varepsilon} + \frac{1}{e}\right) + \tfrac{2}{a_2}((4 \gamma_2 M_h + 1)\tfrac{b_h}{\gamma_1} + \tfrac{b_\tau}{\delta} + 1)\int_\Omega \left(c_{2 \varepsilon} \ln c_{2 \varepsilon} + \frac{1}{e}\right) + \tfrac{b_h}{2 \gamma_1} \int_\Omega \dfrac{|\nabla h_\varepsilon|^2}{h_\varepsilon} + \tfrac{b_\tau}{2 \delta} \int_\Omega \dfrac{|\nabla \tau_\varepsilon|^2}{\tau_\varepsilon}.
\end{equation*}
and
\begin{align*}
  \mathcal{D}_\varepsilon(t) := & \frac{a_1}{4}\int_\Omega \dfrac{|\nabla c_{1 \varepsilon}|^2}{c_{1\varepsilon}} + \frac{1}{2} \int_\Omega \dfrac{|\nabla c_{2 \varepsilon}|^2}{c_{2\varepsilon}} + \dfrac{\varepsilon}{2} \int_\Omega c_{1 \varepsilon}^\theta \ln(2 + c_{1 \varepsilon}) + \dfrac{\beta}{2} \int_\Omega c_{1 \varepsilon}^2 \ln(2 + c_{1 \varepsilon}) \\[5pt]
  & + \frac{b_h}{2}\int_\Omega \dfrac{|\nabla h_\varepsilon|^2}{h_\varepsilon} c_{1\varepsilon} + \frac{b_\tau}{2} \int_\Omega \dfrac{|\nabla \tau_\varepsilon|^2}{\tau_\varepsilon} c_{1\varepsilon}.
\end{align*}
Then, from \eqref{entro1} with the help of \eqref{pent1}, \eqref{pent6} and setting $\eta \leq \min\{\varrho_1,\varrho_2,\mu,\sigma\}$ we have
\begin{equation}
 \label{entro5}
 \mathcal{F}'_\varepsilon(t) + \eta \mathcal{F}_\varepsilon(t) + \mathcal{D}_\varepsilon(t) \leq C_{16} \quad ~~\text{for all} ~  t \in (0, T)
\end{equation}
as $ \mathcal{D}_\varepsilon(t)$ is nonnegative, from \eqref{entro5} we can have the following ordinary differential equation
\begin{equation}
 \label{entro6}
 \mathcal{F}'_\varepsilon(t) + \eta \mathcal{F}_\varepsilon(t)  \leq C_{16} \quad ~~\text{for all} ~  t \in (0, T).
\end{equation}
Lemma 3.4 of \cite{stinner2014global} allow us to find a $C_{14}(T) > 0$ such that
\begin{equation}
  \label{entro7}
   \mathcal{F}_\varepsilon(t)  \leq C_{17}(T) \quad ~~\text{for all} ~  t \in (0, T)
\end{equation}
thus yielding \eqref{entropy_main1}.  A straightforward integration of \eqref{entro5} over $(0, T)$ in view of \eqref{entro7} results in
  \begin{equation*}
    \int_0^T \mathcal{D}_\varepsilon(t) \leq \widehat{C}(T),
  \end{equation*}
   which gives \eqref{entropy_main2}.

  \end{proof}

%%%%%%%%%%%%%%%%%%%%%%%%%%%%%%%%%%%%%%%%%%%%%%%%%%%%%%%%%%%%%%%%%%%%%%%%%%%%%%%%%%%%%%%%%%%%%%%%%%%%%%%%%%%%%%%%%%%%%%%%%%%%%%%%%%%%%%%%%%%%%%%%%%%%%%%%%%%%%%%%%%%%%%%%%%%%%%%%%%%%%%%%%%%%%%%%%%%%%%%%%%%%%%%
%%%%%%%%%%%%%%%%%%%%%%%%%%%%%%%%%%%%%%%%%%%%%%%%%%%%%%%%%%%%%%%%%%%%%%%%%%%%%%%%%%%%%%%%%%%%%%%%%%%%%%%%%%%%%%%%%%%%%%%%%%%%%%%%%%%%%%%%%%%%%%%%%%%%%%%%%%%%%%%%%%%%%%%%%%%%%%%%%%%%%%%%%%%%%%%%%%%%%%%%%%%%%%%
%%%%%%%%%%%%%%%%%%%%%%%%%%%%%%%%%%%%%%%%%%%%%%%%%%%%%%%%%%%%%%%%%%%%%%%%%%%%%%%%%%%%%%%%%%%%%%%%%%%%%%%%%%%%%%%%%%%%%%%%%%%%%%%%%%%%%%%%%%%%%%%%%%%%%%%%%%%%%%%%%%%%%%%%%%%%%%%%%%%%%%%%%%%%%%%%%%%%%%%%%%%%%%%
\subsection*{Further $\varepsilon$-independent estimates}
Based on the inequalities \eqref{entropy_main1} and \eqref{entropy_main2} of Lemma \ref{entropy_lemma13} we will proceed to establish some estimates for the solution components which shall result in strong compactness properties.

 \begin{Lemma}
 \label{nab_bound_lem}
   Let $T > 0$. Then there exists a $C(T) > 0$ such that for every $\varepsilon \in (0, 1)$
   \begin{align}
     \label{h_tau_nab} \int_\Omega |\nabla h_\varepsilon(\cdot, t)|^2 \leq & C(T),  \quad  \int_\Omega |\nabla \tau_\varepsilon(\cdot, t)|^2 \leq C(T),  \quad \text{for all} ~ t \in (0, T),\\[5pt]
     \label{c2_sq} & \int_0^T \int_\Omega c_{2\varepsilon}^2 \leq C(T).
   \end{align}
 \end{Lemma}
\begin{proof}
  To validate the estimates in \eqref{h_tau_nab}, we use \eqref{h_bound}, \eqref{t_bound} and \eqref{entropy_main1} to see that there exists a $C(T) > 0$ such that %{\cmg The 2nd inequalities below only hold when $c_{1\eps}, c_{2\eps}\ge 1$, what if $c_{1\eps}<1$ and/or  $c_{2\eps}<1$?}
  \begin{align*}
    & \int_\Omega |\nabla h_\varepsilon|^2 =  \int_\Omega \frac{|\nabla h_\varepsilon|^2}{h_\varepsilon} \cdot h_\varepsilon \leq M_h \int_\Omega \frac{|\nabla h_\varepsilon|^2}{h_\varepsilon} \leq C(T),  \quad \text{for all} ~ t \in (0, T)    \\[5pt]
    & \int_\Omega |\nabla \tau_\varepsilon|^2 =  \int_\Omega \frac{|\nabla \tau_\varepsilon|^2}{\tau_\varepsilon} \cdot \tau_\varepsilon \leq M_\tau \int_\Omega \frac{|\nabla \tau_\varepsilon|^2}{\tau_\varepsilon} \leq C(T) \quad \text{for all} ~ t \in (0, T).
  \end{align*}
 To prove \eqref{c2_sq}, using the Gagliardo-Nirenberg inequality in conjunction with \eqref{L1bound} we can directly have
 \begin{align}
 \label{c2_sq1}
    \nonumber & \int_\Omega c_{2\varepsilon}^2 \leq C_1 \cdot \left\{\int_\Omega \dfrac{|\nabla c_{2\varepsilon}|^2}{c_{2\varepsilon}}\right\} \cdot \left\{\int_\Omega c_{2\varepsilon}\right\} + C_2 \left\{\int_\Omega c_{2\varepsilon}\right\}^2  \\[5pt]
    & \int_\Omega c_{2\varepsilon}^2 \leq C_3(T) \int_\Omega \dfrac{|\nabla c_{2\varepsilon}|^2}{c_{2\varepsilon}} + C_4(T)  \quad \text{for all} ~  t \in (0, T).
  \end{align}
 An integration of \eqref{c2_sq1} in view of \eqref{entropy_main2} results in
  \begin{align}
 \label{c2_sq2}
    \int_0^T \int_\Omega c_{2\varepsilon}^2 \leq C_3(T) \int_0^T \int_\Omega \dfrac{|\nabla c_{2\varepsilon}|^2}{c_{2\varepsilon}} + C_4(T) T \leq C(T).
  \end{align}
\end{proof}

%%%%%%%%%%%%%%%%%%%%%%%%%%%%%%%%%%%%%%%%%%%%%%%%%%%%%%%%%%%%%%%%%%%%%%%%%%%%%%%%%%%%%%%%%%%%%%%%%%%%%%%%%%%%%%%%%%%%%%%%%%%%%%%%%%%%%%%%%%%%%%%%%%%%%%%%%%%%%%%%%%%%%%%%%%%%%%%%%%%%%%%%%%%%%%%%%%%%%%%%%%%%%%%
\begin{Lemma}
  \label{eind_ht}
  Let $T > 0$ and assume that $k > \frac{n + 2}{2}$ is an integer, then there exists a $C(T) > 0$ such that for all $\varepsilon \in (0, 1)$
\begin{align}
    & \label{h_com_bound} \|\partial_t h_\varepsilon\|_{L^1((0, T); (W_0^{k, 2}(\Omega))^\ast)} \leq C(T),                \\[5pt]
    & \label{tau_com_bound} \|\partial_t \tau_\varepsilon\|_{L^1((0, T); (W_0^{k, 2}(\Omega))^\ast)} \leq C(T).
  \end{align}
\end{Lemma}
\begin{proof}
 To derive \eqref{h_com_bound}, by utilizing  \eqref{L1bound}, \eqref{h_bound}, \eqref{entropy_main1} and the fact that $0 \leq \tfrac{c_{2\varepsilon}}{1 + c_{2\varepsilon}} \leq 1$ for $c_{2\varepsilon} \geq 0$, we can find a $C_1(T) > 0$ such that
 \begin{align}
  \label{hpre2}
    \nonumber \left|\int_\Omega \partial_t h_\varepsilon \cdot \phi\right| & = \left| \varepsilon \int_\Omega \nabla h_\varepsilon \cdot \nabla \phi - \gamma_1\int_\Omega h_\varepsilon c_{1\varepsilon}\phi  - \gamma_2\int_\Omega h_\varepsilon c_{2\varepsilon}\phi - \mu \int_\Omega h_\varepsilon \phi +  \int_\Omega \dfrac{c_{2\varepsilon}}{1 + c_{2\varepsilon}} \phi\right|  \\[5pt]
    & \nonumber \leq \varepsilon \left(\int_\Omega \dfrac{|\nabla h_\varepsilon|^2}{h_\varepsilon}\right)^{\frac{1}{2}} \cdot \left(\int_\Omega h_\varepsilon\right)^{\frac{1}{2}} \|\nabla \phi\|_{L^\infty(\Omega)}
    + \gamma_1 \left(\int_\Omega h_\varepsilon c_{1\varepsilon}\right)\|\phi\|_{L^\infty(\Omega)}\\[7pt]
    & \nonumber  + \gamma_2 \left(\int_\Omega h_\varepsilon c_{2\varepsilon}\right)\|\phi\|_{L^\infty(\Omega)} + \mu \left(\int_\Omega h_\varepsilon \right)\|\phi\|_{L^\infty(\Omega)} + \left(\int_\Omega\dfrac{c_{2\varepsilon}}{1 + c_{2\varepsilon}}\right) \|\phi\|_{L^\infty(\Omega)} \\[7pt]
    & \nonumber \leq \varepsilon \dfrac{\sqrt{M_h |\Omega|}}{2} \left(\int_\Omega \dfrac{|\nabla h_\varepsilon|^2}{h_\varepsilon} + 1\right)  \|\nabla \phi\|_{L^\infty(\Omega)} + \gamma_1 M_h \left(\int_\Omega c_{1\varepsilon}\right)\|\phi\|_{L^\infty(\Omega)} \\[7pt]
    & \nonumber + \gamma_2 M_h \left(\int_\Omega c_{2\varepsilon}\right)\|\phi\|_{L^\infty(\Omega)} + \mu M_h|\Omega| \|\phi\|_{L^\infty(\Omega)} + |\Omega| \|\phi\|_{L^\infty(\Omega)} \\[7pt]
    & \leq C_1\|\phi\|_{W^{1, \infty}(\Omega)} ~~\text{for all}~ t \in (0, T) ~ \text{and each}~ \phi \in C^\infty_{0}(\Omega).
  \end{align}
 The Sobolev embedding theorem allows us to have
  \begin{equation}
   \label{hpre3}
    \|\phi\|_{W^{1, \infty}(\Omega)} \leq C_2\|\phi\|_{W^{k, 2}_0(\Omega)} ~~\text{for all}~ \phi \in C^\infty_{0}(\Omega)
  \end{equation}
  with $C_2 >0$. From \eqref{hpre2} and \eqref{hpre3} we can directly have
  \begin{equation*}
   \|\partial_t h_\varepsilon\|_{L^1((0, T); (W_0^{k, 2}(\Omega))^\ast)} = \int_0^T \sup_{\substack{\phi \in C_0^\infty(\Omega),\\ \|\phi\|_{W^{k, 2}_0(\Omega)} = 1}} \left|\int_\Omega \partial_t h_\varepsilon \cdot \phi\right| dt \leq C_1 C_2 T
  \end{equation*}
  thus establishing \eqref{h_com_bound}. Analogously, we can derive \eqref{tau_com_bound}.

\end{proof}

%%%%%%%%%%%%%%%%%%%%%%%%%%%%%%%%%%%%%%%%%%%%%%%%%%%%%%%%%%%%%%%%%%%%%%%%%%%%%%%%%%%%%%%%%%%%%%%%%%%%%%%%%%%%%%%%%%%%%%%%%%%%%%%%%%%%%%%%%%%%%%%%%%%%%%%%%%%%%%%%%%%%%%%%%%%%%%%%%%%%%%%%%%%%%%%%%%%%%%%%%%%%%%%
\begin{Lemma}
  \label{c1_tim_der}
  Let $T > 0$. Suppose $l > \frac{n + 4}{2}$ is an integer then there exists a $C(T) > 0$ such that for all $\varepsilon \in (0, 1)$
  \begin{align}
    \label{c1_tim1} & \|c_{1\varepsilon}\|_{L^{\frac{4}{3}}((0, T); W^{1, \frac{4}{3}}(\Omega))} \leq C(T),   \\[5pt]
    \label{c1_tim2} & \|\partial_t c_{1\varepsilon}\|_{L^{1}((0, T); (W^{l, 2}_0(\Omega))^\ast)} \leq C(T), \quad \text{and}, \\[5pt]
    \label{c2_tim1} & \|c_{2\varepsilon}\|_{L^{\frac{5}{4}}((0, T); W^{1, \frac{5}{4}}(\Omega))} \leq C(T),   \\[5pt]
    \label{c2_tim2} & \|\partial_t c_{2\varepsilon}\|_{L^{1}((0, T); (W^{l, 2}_0(\Omega))^\ast)} \leq C(T).
  \end{align}
\end{Lemma}

\begin{proof}
  From \eqref{sqL1bound} and \eqref{entropy_main2} we infer that there exists a $C_1(T) > 0$ fulfilling
  \begin{align*}
    \int_0^{T} \int_\Omega c_{1\varepsilon}^2 \leq C_1(T) \quad ~\text{and}~ \quad \int_0^T \int_\Omega \dfrac{|\nabla c_{1 \varepsilon}|^2}{c_{1\varepsilon}} \leq C_1(T)
  \end{align*}
  this directly yields with the H\"older inequality
  \begin{equation}
    \label{c1_tres1}
    \int_0^T \int_\Omega |\nabla c_{1\varepsilon}|^{\frac{4}{3}} \leq \left(\int_0^T \int_\Omega \dfrac{|\nabla c_{1 \varepsilon}|^2}{c_{1\varepsilon}}\right)^{\frac{2}{3}} \cdot  \left( \int_0^{T} \int_\Omega c_{1\varepsilon}^2 \right)^{\frac{1}{3}} \leq C_2(T),
  \end{equation}
 which ensures that  $\|c_{1\varepsilon}\|_{L^{\frac{4}{3}}(\Omega \times (0, T))}$ is bounded, thus establishing \eqref{c1_tim1}.\\[-2ex]

\noindent
To prove \eqref{c1_tim2}, by using \eqref{alpha_def}, \eqref{app2}, \eqref{t_bound} and Young's inequality we can have for all $\phi \in C_0^\infty(\Omega)$
  \begin{align}
  \label{c1_tim5}
    \nonumber \left|\int_\Omega \partial_t c_{1\varepsilon} \cdot \phi\right| & = \bigg| - a_1 \int_\Omega \nabla c_{1\varepsilon} \cdot \nabla \phi + b_h \int_\Omega c_{1\varepsilon} \nabla h_\varepsilon \cdot \nabla \phi + b_\tau \int_\Omega c_{1\varepsilon} \nabla \tau_\varepsilon \cdot \nabla \phi \\[5pt]
    \nonumber &  - \int_\Omega \alpha_1(\tau_\varepsilon)c_{1\varepsilon} \phi  + \int_\Omega \alpha_2( \tau_\varepsilon) F_\varepsilon(c_{2\varepsilon})\phi + \beta \int_\Omega  c_{1\varepsilon}\phi\\[5pt]
    \nonumber &  -  \beta \int_\Omega  c_{1\varepsilon}^2\phi - \beta \int_\Omega  c_{1\varepsilon}   c_{2\varepsilon}\phi  - \beta \int_\Omega  c_{1\varepsilon} \tau_\varepsilon\phi  - \varepsilon \int_\Omega  c_{1\varepsilon}^\theta \phi \bigg|  \\[5pt]
    \nonumber & \leq a_1 \left(\frac{3}{4} \int_\Omega |\nabla c_{1\varepsilon}|^{\frac{4}{3}} + \frac{1}{4} |\Omega|\right)\|\nabla \phi\|_{L^\infty(\Omega)} + \frac{b_h}{2} \left(\int_\Omega c_{1\varepsilon}^2 + \int_\Omega |\nabla h_\varepsilon|^2\right)\|\nabla \phi\|_{L^\infty(\Omega)} \\[5pt]
    \nonumber & + \frac{b_\tau}{2} \left(\int_\Omega c_{1\varepsilon}^2 + \int_\Omega |\nabla \tau_\varepsilon|^2\right)\|\nabla \phi\|_{L^\infty(\Omega)} + M_{\alpha_1} \left(\int_\Omega c_{1\varepsilon}\right)\| \phi\|_{L^\infty(\Omega)} \\[5pt]
    \nonumber & +  M_{\alpha_2} \left(\int_\Omega c_{2\varepsilon}\right)\| \phi\|_{L^\infty(\Omega)} + \beta \left(\int_\Omega c_{1\varepsilon}\right)\| \phi\|_{L^\infty(\Omega)} + \beta \left(\int_\Omega c_{1\varepsilon}^2\right)\| \phi\|_{L^\infty(\Omega)} \\[5pt]
    \nonumber & + \frac{\beta}{2} \left(\int_\Omega c_{1\varepsilon}^2 + \int_\Omega c_{2\varepsilon}^2\right)\| \phi\|_{L^\infty(\Omega)} + \beta M_\tau \left(\int_\Omega c_{1\varepsilon}\right)\| \phi\|_{L^\infty(\Omega)} + \varepsilon \left(\int_\Omega c_{1\varepsilon}^\theta\right)\| \phi\|_{L^\infty(\Omega)} \\[5pt]
    & \leq Z_\varepsilon(t) \|\phi\|_{W^{2, \infty}(\Omega)} \quad \text{for all} ~ t \in (0, T),
  \end{align}
  where
  \begin{align*}
    Z_\varepsilon(t) := & \frac{3}{4} a_1 \int_\Omega |\nabla c_{1\varepsilon}|^{\frac{4}{3}} + \frac{b_h}{2} \int_\Omega |\nabla h_\varepsilon|^{2} + \frac{b_\tau}{2} \int_\Omega |\nabla \tau_\varepsilon|^{2} + \left(\frac{b_h}{2} + \frac{b_\tau}{2} + \frac{3\beta}{2}\right) \int_\Omega c_{1\varepsilon}^2    \\[5pt]
    & + \frac{\beta}{2} \int_\Omega c_{2\varepsilon}^2  + (M_{\alpha_1} + \beta M_\tau + \beta) \int_\Omega c_{1\varepsilon} + M_{\alpha_2} \int_\Omega c_{2\varepsilon} + \varepsilon \int_\Omega c_{1\varepsilon}^\theta + \frac{1}{4}a_1|\Omega|
  \end{align*}
for all $\varepsilon \in (0, 1)$ and each $t \in (0, T)$. Using \eqref{L1bound}, \eqref{sqL1bound}, \eqref{epL1bound}, \eqref{h_tau_nab}, \eqref{c2_sq} and \eqref{c1_tres1}, we can find a $C_3(T) > 0$ such that
\begin{equation}
 \label{c1_tim6}
  \int_0^T Z_\varepsilon({\cb t}) dt \leq C_3(T).
\end{equation}
\noindent
From \eqref{c1_tim5}, \eqref{c1_tim6}, and the Sobolev embedding $W_0^{l, 2}(\Omega) \hookrightarrow W^{2, \infty}(\Omega)$ we can have
\begin{align*}
   \|\partial_t c_{1\varepsilon}\|_{L^1((0, T); (W_0^{l, 2}(\Omega))^\ast)} & = \int_0^T \sup_{\substack{\phi \in C_0^\infty(\Omega),\\ \|\phi\|_{W^{l, 2}_0(\Omega)} = 1}} \left|\int_\Omega \partial_t c_{1\varepsilon} \cdot \phi \right| dt  \leq \int_0^T \sup_{\substack{\phi \in C_0^\infty(\Omega),\\ \|\phi\|_{W^{l, 2}_0(\Omega)} = 1}} Z_\varepsilon(t) \|\phi\|_{W^{2, \infty}(\Omega)} dt  \\[5pt]
   & \leq C_3(T) C_4,
  \end{align*}
which entails \eqref{c1_tim2}\\[-2ex]

\noindent
Integrating \eqref{ent5} and using \eqref{entropy_main2} we can have
%\begin{align}
%\label{c2_tim6}
%  \nonumber\int_\Omega c_{2\varepsilon}^{\frac{5}{3}} & \leq C_5 \|\nabla \sqrt{ c_{2\varepsilon}}\|_{L^2(\Omega)}^{\frac{2n}{3}} \|\sqrt{c_{2\varepsilon}}\|^{\frac{10 - 2n}{3}}_{L^2(\Omega)} + C_6 \|\sqrt{c_{2\varepsilon}}\|^{\frac{10}{3}}_{L^2(\Omega)}                   \\[5pt]
%  \nonumber  & \leq C_7 \left(\int_\Omega |\nabla \sqrt{ c_{2\varepsilon}}|^2 \right)^{\frac{n}{3}} + C_7   \\[5pt]
%  & \leq C_8 \int_\Omega |\nabla \sqrt{ c_{2\varepsilon}}|^2 + C_8 \leq C_9 \int_\Omega \frac{|\nabla c_{2\varepsilon}|^2}{c_{2\varepsilon}} + C_9 \quad \text{for all} ~ t \in (0, T).
%\end{align}
%Integrating \eqref{c2_tim6}, keeping in mind \eqref{entropy_main2} we can have
\begin{align}
  \label{c2_tim8}
  \int_0^T \int_\Omega c_{2\varepsilon}^{\frac{5}{3}}  \leq 
  %C_9 \int_0^T \int_\Omega \frac{|\nabla c_{2\varepsilon}|^2}{c_{2\varepsilon}} + 
  C_5(T).
\end{align}
\eqref{c2_tim8} taken together with \eqref{entropy_main2}, along with H\"older's inequality, allow us to have
\begin{equation}
\label{c2_tim9}
  \int_0^T \int_\Omega |\nabla c_{2\varepsilon}|^{\frac{5}{4}} \leq \left(\int_0^T \int_\Omega \frac{|\nabla c_{2\varepsilon}|^2}{c_{2\varepsilon}}\right)^{\frac{5}{8}} \cdot \left(\int_0^T \int_\Omega c_{2\varepsilon}^{\frac{5}{3}}\right)^{\frac{3}{8}} \leq C_{6}(T),
\end{equation}
%By H\"{o}lder's inequality 
thus we can deduce \eqref{c2_tim1} from \eqref{c2_tim8} and \eqref{c2_tim9}.\\[-2ex]

\noindent
To prove \eqref{c2_tim2}, we use \eqref{alpha_def}, \eqref{app2} and Young's inequality to get for all $\phi \in C_0^\infty(\Omega)$
  \begin{align}
  \label{c2_tim5}
    \nonumber \left|\int_\Omega \partial_t c_{2\varepsilon} \cdot \phi\right| & = \bigg| - a_2 \int_\Omega \nabla c_{2\varepsilon} \cdot \nabla \phi  + \int_\Omega \alpha_1(\tau_\varepsilon)c_{1\varepsilon} \phi - \int_\Omega \alpha_2( \tau_\varepsilon) F_\varepsilon(c_{2\varepsilon}) \phi  \bigg| \\[5pt]
    \nonumber & \leq a_2 \left(\frac{4}{5} \int_\Omega |\nabla c_{2\varepsilon}|^{\frac{5}{4}} + \frac{1}{5} |\Omega|\right)\|\nabla \phi\|_{L^\infty(\Omega)} + M_{\alpha_1} \left(\int_\Omega c_{1\varepsilon}\right)\| \phi\|_{L^\infty(\Omega)} \\[5pt]
    \nonumber & +  M_{\alpha_2} \left(\int_\Omega c_{2\varepsilon}\right)\| \phi\|_{L^\infty(\Omega)}\\[5pt]
    & \leq V_\varepsilon(t) \|\phi\|_{W^{2, \infty}(\Omega)} \quad \text{for all} ~ t \in (0, T),
  \end{align}
   where,
  \begin{align*}
    V_\varepsilon(t) := & \frac{4}{5} a_2 \int_\Omega |\nabla c_{2\varepsilon}|^{\frac{5}{4}} + M_{\alpha_1} \int_\Omega c_{1\varepsilon} + M_{\alpha_2} \int_\Omega c_{2\varepsilon} + \frac{1}{5} a_2 |\Omega|
  \end{align*}
for all $\varepsilon \in (0, 1)$ and each $t \in (0, T)$. Using \eqref{L1bound} and \eqref{c2_tim9}, we can find a $C_{7}(T) > 0$ such that
\begin{equation}
 \label{c2_tim611}
  \int_0^T V_\varepsilon({\cb t}) dt \leq C_{7}(T).
\end{equation}
From \eqref{c2_tim5}, \eqref{c2_tim611}, and the Sobolev embedding $W_0^{l, 2}(\Omega) \hookrightarrow W^{2, \infty}(\Omega)$ we can then have
\begin{align*}
   \|\partial_t c_{2\varepsilon}\|_{L^1((0, T); (W_0^{l, 2}(\Omega))^\ast)} & = \int_0^T \sup_{\substack{\phi \in C_0^\infty(\Omega),\\ \|\phi\|_{W^{l, 2}_0(\Omega)} = 1}} \left|\int_\Omega \partial_t c_{2\varepsilon} \cdot \phi \right| dt  \leq \int_0^T \sup_{\substack{\phi \in C_0^\infty(\Omega),\\ \|\phi\|_{W^{l, 2}_0(\Omega)} = 1}} V_\varepsilon(t) \|\phi\|_{W^{2, \infty}(\Omega)} dt  \\[5pt]
   & \leq C_{7}(T) C_{8},
  \end{align*}
which entails \eqref{c2_tim2}.

\end{proof}
%%%%%%%%%%%%%%%%%%%%%%%%%%%%%%%%%%%%%%%%%%%%%%%%%%%%%%%%%%%%%%%%%%%%%%%%%%%%%%%%%%%%%%%%%%%%%%%%%%%%%%%%%%%%%%%%%%%%%%%%%%%%%%%%%%%%%%%%%%%%%%%%%%%%%%%%%%%%%%%%%%%%%%%%%%%%%%%%%%%%%%%%%%%%%%%%%%%%%%%%%%%%%%%%%%%%

\noindent
Using Lemma \ref{eind_ht} and Lemma \ref{c1_tim_der} in conjunction with the Aubin-Lions lemma \cite[Chapter 3]{temam2001navier}, we can establish the following strong precompactness properties.

\begin{Lemma}
\label{precompactness_all}
Let $T > 0$. Then
\begin{align}
  \label{precompact_c1}    & \{c_{1\varepsilon}\}_{\varepsilon \in (0,1)} ~\text{is strongly precompact in}~ L^{\frac{4}{3}}(\Omega \times (0, T)),  \\[5pt]
  \label{precompact_c2}    & \{c_{2\varepsilon}\}_{\varepsilon \in (0,1)} ~\text{is strongly precompact in}~ L^{\frac{5}{4}}(\Omega \times (0, T)),  \\[5pt]
  \label{precompact_h}     & \{h_{\varepsilon}\}_{\varepsilon \in (0,1)}   ~\text{is strongly precompact in}~ L^{2}(\Omega \times (0, T)),   \quad ~\text{and},        \\[5pt]
  \label{precompact_tau}   & \{\tau_{\varepsilon}\}_{\varepsilon \in (0,1)}   ~\text{is strongly precompact in}~ L^{2}(\Omega \times (0, T)).
\end{align}
\end{Lemma}

\begin{proof}
  Let $l$ be the integer chosen in \eqref{c1_tim2}. As  $(W_0^{l, 2}(\Omega))^\ast$ is a Hilbert space and the embedding $W^{1, \frac{4}{3}}(\Omega) \hookrightarrow L^{\frac{4}{3}}(\Omega)$ is compact, we can invoke the Aubin-Lions lemma \cite[Theorem 3.2.3 and Remark 3.2.1]{temam2001navier}, \eqref{c1_tim1} and \eqref{c1_tim2} to conclude that $\{c_{1\varepsilon}\}_{\varepsilon \in (0,1)}$ is strongly precompact in $L^{\frac{4}{3}}(\Omega \times (0, T))$ thus establishing \eqref{precompact_c1}. We can verify \eqref{precompact_c2}, by collecting \eqref{c2_tim1} and \eqref{c2_tim2} and again invoking the Aubin-Lions lemma. Similarly, \eqref{precompact_h} and \eqref{precompact_tau} can be established by again using the Aubin-Lions lemma and collecting \eqref{h_tau_nab}, \eqref{h_com_bound}-\eqref{tau_com_bound}.

\end{proof}

%%%%%%%%%%%%%%%%%%%%%%%%%%%%%%%%%%%%%%%%%%%%%%%%%%%%%%%%%%%%%%%%%%%%%%%%%%%%%%%%%%%%%%%%%%%%%%%%%%%%%%%%%%%%%%%%%%%%%%%%%%%%%%%%%%%%%%%%%%%%%%%%%%%%%%%%%%%%%%%%%%%%%%%%%%%%%%%%%%%%%%%%%%%%%%%%%%%%%%%%%%%%%%%%%%%%

\section{Construction of weak solutions}

\begin{proof}[Proof of Theorem \ref{main_theorem}]
Estimate \eqref{entropy_main2} allows us to claim the boundedness of $\{c_{1\varepsilon}^2 \ln (2 + c_{1\varepsilon})\}_{\varepsilon \in (0, 1)}$ and $\{\varepsilon c_{1\varepsilon}^\theta \ln (2 + c_{1\varepsilon})\}_{\varepsilon \in (0, 1)}$ in the space $L^1_{loc}(\bar{\Omega} \times [0, \infty))$, and hence, $\{c_{1\varepsilon}^2\}_{\varepsilon \in (0, 1)}$ and  $\{\varepsilon  c_{1\varepsilon}^\theta\}_{\varepsilon \in (0, 1)}$ are equi-integrable. By the Dunford-Pettis theorem \cite[A8.14]{alt2016linear} we can claim that $\{c_{1\varepsilon}^2\}_{\varepsilon \in (0, 1)}$ and  $\{\varepsilon  c_{1\varepsilon}^\theta\}_{\varepsilon \in (0, 1)}$ are weakly sequentially precompact in $L^1_{loc}(\bar{\Omega} \times [0, \infty))$. Also, from \eqref{entropy_main2} we can infer that $\{c_{1\varepsilon}\}_{\varepsilon \in (0, 1)}$ is bounded in $L^2_{loc}(\bar{\Omega} \times [0, \infty))$. These in conjunction with \eqref{c1_tim1} and \eqref{precompact_c1} allow us to apply a standard extraction procedure to select a suitable subsequence $\{\varepsilon_j\}_{j \in \mathbb{N}}$ such that
\begin{equation*}
  \varepsilon_j \in (0, 1) ~\text{for all}~ j \in \mathbb{N}, ~\text{and}~ \varepsilon_j \searrow 0 ~\text{as}~ j \to \infty
\end{equation*}
and a nonnegative function
\begin{equation*}
\label{th3c1_main}
c_1 \in L^2_{loc}(\bar{\Omega} \times [0, \infty)) \cap L_{loc}^{\frac{4}{3}}([0, \infty) ; W^{1, \frac{4}{3}}(\Omega))
\end{equation*}
\begin{equation}
\begin{cases}
 \label{mth1}
  c_{1\varepsilon_j} & \to c_1 ~\text{a.e. in}~ \Omega \times (0, \infty),                                                                          \\[5pt]
  c_{1\varepsilon_j}^2 & \rightharpoonup c_1^2 ~\text{in}~ L^1_{loc}(\bar{\Omega} \times [0, \infty)),                                              \\[5pt]
  c_{1\varepsilon_j} & \rightharpoonup c_1 ~\text{in}~ L^2_{loc}(\bar{\Omega} \times [0, \infty)),                                                  \\[5pt]
  \nabla c_{1\varepsilon_j} & \rightharpoonup \nabla c_1 ~\text{in}~ L^{\frac{4}{3}}_{loc}(\bar{\Omega} \times [0, \infty)),                        \\[5pt]
  \varepsilon_j c_{1\varepsilon_j}^\theta & \to 0 ~\text{a.e. in}~ \Omega \times (0, \infty),      \quad \text{and}                                    \\[5pt]
  \varepsilon_j c_{1\varepsilon_j}^\theta & \rightharpoonup 0 ~\text{in}~ L_{loc}^1(\bar{\Omega} \times [0, \infty))
\end{cases}
\end{equation}
as $j \to \infty$. With $1$ as test function, from \eqref{mth1} we can directly have
\begin{equation}
  \label{mth411}
   c_{1\varepsilon_j} \to c_1 ~\text{in}~ L^2_{loc}(\bar{\Omega} \times [0, \infty)) ~ \text{as} ~ j \to \infty.
\end{equation}
Similarly, from \eqref{entropy_main1}, \eqref{entropy_main2}, \eqref{h_tau_nab}, \eqref{c2_tim1}, \eqref{precompact_c2}-\eqref{precompact_tau}, passing to a subsequence if necessary, there exist some nonnegative functions $c_2$, $h$, and $\tau$ with
\begin{equation}
	\label{2mth1}
	\begin{cases}
		& c_2 \in L^{\frac{5}{4}}_{loc}([0, \infty); W^{1, \frac{5}{4}}(\Omega)),                          \\[5pt]
		& h \in L^\infty(\Omega \times (0, \infty)) \cap L^2_{loc}([0, \infty) ; W^{1, 2}(\Omega)),        \\[5pt]
		& \tau \in L^\infty(\Omega \times (0, \infty)) \cap L^2_{loc}([0, \infty) ; W^{1, 2}(\Omega)).
	\end{cases}
\end{equation}
such that
\begin{align}
  \label{mth4} & c_{2\varepsilon_j} \to c_2 ~\text{in}~ L^{\frac{5}{4}}_{loc}(\bar{\Omega} \times [0, \infty)) ~\text{and a.e. in}~ \Omega \times (0, \infty),                        \\[5pt]
  \label{mth5} & \nabla c_{2\varepsilon_j}  \rightharpoonup \nabla c_2 ~\text{in}~ L^{\frac{5}{4}}_{loc}(\bar{\Omega} \times [0, \infty)),                                            \\[5pt]
  \label{mth6} & h_{\varepsilon_j} \to h ~\text{in}~ L^{2}_{loc}(\bar{\Omega} \times [0, \infty)) ~\text{and a.e. in}~ \Omega \times (0, \infty),                                     \\[5pt]
  \label{mth7} & \nabla h_{\varepsilon_j}  \rightharpoonup \nabla h ~\text{in}~ L^{2}_{loc}(\bar{\Omega} \times [0, \infty)),                                                         \\[5pt]
  \label{mth8} & \tau_{\varepsilon_j} \to \tau ~\text{in}~ L^{2}_{loc}(\bar{\Omega} \times [0, \infty)) ~\text{and a.e. in}~ \Omega \times (0, \infty),  \quad ~\text{and}            \\[5pt]
  \label{mth9} & \nabla \tau_{\varepsilon_j}  \rightharpoonup \nabla \tau ~\text{in}~ L^{2}_{loc}(\bar{\Omega} \times [0, \infty))
\end{align}
as $j \to \infty$.

\noindent
Let $T > 0$ and $\phi \in C^\infty_0(\bar{\Omega} \times [0, T))$ with $\dfrac{\partial \phi}{\partial \nu} = 0$ on $\partial \Omega \times [0, T)$. Then from the first equation of \eqref{model1} we have for all $\varepsilon \in (0, 1)$

\begin{align}
   \label{mthdefeq1}
  \nonumber  & - \int_0^T \int_\Omega c_{1\varepsilon} \partial_t \phi - \int_\Omega c_{10\varepsilon} \phi(\cdot, 0) = - a_1 \int_0^T \int_\Omega \nabla c_{1\varepsilon} \cdot \nabla \phi + b_h \int_0^T \int_\Omega c_{1\varepsilon} \nabla h_\varepsilon \cdot \nabla \phi \\[5pt]
  \nonumber  & + b_\tau \int_0^T \int_\Omega c_{1\varepsilon} \nabla \tau_\varepsilon \cdot \nabla \phi - \int_0^T \int_\Omega \alpha_1(\tau_\varepsilon )c_{1\varepsilon} \phi + \int_0^T \int_\Omega \alpha_2(\tau_\varepsilon )F_\varepsilon(c_{2\varepsilon}) \phi  + \beta \int_0^T \int_\Omega c_{1\varepsilon} \phi\\[5pt]
  & - \beta \int_0^T \int_\Omega  c_{1\varepsilon}^2 \phi - \beta \int_0^T \int_\Omega  c_{1\varepsilon} c_{2\varepsilon} \phi - \beta \int_0^T \int_\Omega \tau_\varepsilon c_{1\varepsilon}\phi - \varepsilon \int_0^T \int_\Omega c_{1\varepsilon}^\theta \phi.
\end{align}
By \eqref{reg_assumptions}, \eqref{mth1} and \eqref{mth411} we can deduce that
\begin{equation}
  \label{weak_sol1}
  - \int_0^T \int_\Omega c_{1\varepsilon} \partial_t \phi - \int_\Omega c_{10\varepsilon} \phi(\cdot, 0) \to - \int_0^T \int_\Omega c_{1} \partial_t \phi - \int_\Omega c_{10} \phi(\cdot, 0),
\end{equation}
and
\begin{equation}
\label{weak_sol2}
 - a_1 \int_0^T \int_\Omega \nabla c_{1\varepsilon} \cdot \nabla \phi \to  - a_1 \int_0^T \int_\Omega \nabla c_{1} \cdot \nabla \phi,
\end{equation}
as well as
\begin{align}
\label{weak_sol3}
 - \varepsilon \int_0^T \int_\Omega c_{1\varepsilon}^\theta \phi \to 0
\end{align}
as $\varepsilon = \varepsilon_j \searrow 0$. Combining \eqref{mth1}, \eqref{mth411}, \eqref{mth7} and \eqref{mth9} we can claim that
\begin{align}
\label{weak_sol4}
 b_h \int_0^T \int_\Omega c_{1\varepsilon} \nabla h_\varepsilon \cdot \nabla \phi \to b_h \int_0^T \int_\Omega c_{1} \nabla h \cdot \nabla \phi
\end{align}
and
\begin{align}
 \label{weak_sol5}
 b_\tau \int_0^T \int_\Omega c_{1\varepsilon} \nabla \tau_\varepsilon \cdot \nabla \phi \to b_\tau \int_0^T \int_\Omega c_{1} \nabla \tau \cdot \nabla \phi
\end{align}
as $\varepsilon = \varepsilon_j \searrow 0$. Taken together \eqref{alpha_def}, \eqref{mth1} and \eqref{mth411} gives us
\begin{align}
 \label{weak_sol6}
 - \int_0^T \int_\Omega \alpha_1(\tau_\varepsilon )c_{1\varepsilon} \phi \to - \int_0^T \int_\Omega \alpha_1(\tau) c_{1} \phi,
\end{align}
and \eqref{alpha_def}, \eqref{mth4} with the Vitali convergence theorem allow us to have
\begin{align}
 \label{weak_sol7}
\int_0^T \int_\Omega \alpha_2(\tau_\varepsilon )F_\varepsilon(c_{2\varepsilon}) \phi \to \int_0^T \int_\Omega \alpha_2(\tau )c_{2}\phi.
\end{align}
as $\varepsilon = \varepsilon_j \searrow 0$. Taken together,  \eqref{mth1}, \eqref{mth411} and \eqref{mth4} result in
\begin{align*}
  c_{1\varepsilon} c_{2\varepsilon} \to c_{1} c_{2} \quad \text{in}~L^{\frac{5}{4}}(\Omega \times (0, T)),
\end{align*}
as $\varepsilon = \varepsilon_j \searrow 0$, which implies that
\begin{align}
 \label{weak_sol8}
  - \beta \int_0^T \int_\Omega c_{1\varepsilon}c_{2\varepsilon} \phi \to -\beta \int_0^T \int_\Omega c_{1} c_{2} \phi.
\end{align}
Combining \eqref{t_bound}, \eqref{mth1}, \eqref{mth411} and \cite[Lemma 3.9]{tao2021global}
\begin{align}
 \label{weak_sol9}
  - \beta \int_0^T \int_\Omega c_{1\varepsilon} \tau_\varepsilon \phi \to -\beta \int_0^T \int_\Omega c_{1} \tau \phi.
\end{align}
Finally \eqref{mth1} and \eqref{mth411} taken together yield
\begin{align}
  \label{cweak_sol1}
   \beta \int_0^T \int_\Omega c_{1\varepsilon} \phi - \beta \int_0^T \int_\Omega  c_{1\varepsilon}^2 \phi \to  \beta \int_0^T \int_\Omega c_{1} \phi - \beta \int_0^T \int_\Omega  c_{1}^2 \phi
\end{align}
as $\varepsilon = \varepsilon_j \searrow 0$. The convergence results \eqref{weak_sol1}-\eqref{cweak_sol1} enable us to pass the limit $\varepsilon = \varepsilon_j \searrow 0$ in \eqref{mthdefeq1} to have
\begin{align}
   \label{wmthdefeq1}
  \nonumber  & - \int_0^T \int_\Omega c_{1} \partial_t \phi - \int_\Omega c_{10} \phi(\cdot, 0) = - a_1 \int_0^T \int_\Omega \nabla c_{1} \cdot \nabla \phi + b_h \int_0^T \int_\Omega c_{1} \nabla h \cdot \nabla \phi \\[5pt]
  \nonumber & + b_\tau \int_0^T \int_\Omega c_{1} \nabla \tau \cdot \nabla \phi - \int_0^T \int_\Omega \alpha_1(\tau)c_{1} \phi + \int_0^T \int_\Omega \alpha_2(\tau)c_{2} \phi  \\[5pt]
  & + \beta \int_0^T \int_\Omega c_{1}(1 - c_{1}  -  c_{2}  - \tau )\phi.
\end{align}
Similarly, we can deduce
\begin{align}
   \label{mthdefeq2}
  - \int_0^T \int_\Omega c_{2} \partial_t \phi - \int_\Omega c_{20} \phi(\cdot, 0) = &- a_2 \int_0^T \int_\Omega \nabla c_{2} \cdot \nabla \phi + \int_0^T \int_\Omega \alpha_1(\tau)c_{1} \phi - \int_0^T \int_\Omega \alpha_2(\tau) c_{2} \phi.
\end{align}
\noindent
From the third equation of \eqref{model1} we have for all $\varepsilon \in (0, 1)$
\begin{align}
   \label{mthdefeq3}
  \nonumber  - \int_0^T \int_\Omega h_{\varepsilon} \partial_t \phi - \int_\Omega h_{0\varepsilon} \phi(\cdot, 0) & = - \varepsilon \int_0^T \int_\Omega \nabla h_{\varepsilon} \cdot \nabla \phi
  - \gamma_1 \int_0^T \int_\Omega  h_\varepsilon c_{1\varepsilon} \phi  - \gamma_2 \int_0^T \int_\Omega h_\varepsilon c_{2\varepsilon}\phi  \\[5pt]
  &- \mu \int_0^T \int_\Omega  h_\varepsilon \phi  +  \int_0^T \int_\Omega \dfrac{c_{2\varepsilon}}{1 + c_{2\varepsilon}} \phi.
\end{align}
First we will handle the artificial term on the right-hand side. It follows from \eqref{h_bound}, \eqref{entropy_main1} and the Cauchy-Schwarz inequality that
\begin{align}
\label{artificial1}
  \nonumber \left| - \varepsilon \int_0^T \int_\Omega \nabla h_\varepsilon \cdot \nabla \phi \right| & \leq  \varepsilon \left(\int_0^T \int_\Omega \frac{|\nabla h_\varepsilon|^2}{h_\varepsilon}\right)^{\frac{1}{2}}\left(\int_0^T \int_\Omega h_\varepsilon\right)^{\frac{1}{2}}\|\nabla \phi\|_{L^\infty(\Omega)} \\[5pt]
  \nonumber  & \leq \varepsilon T \sqrt{M_h |\Omega|}\left(\sup_{t \in (0, T)}  \int_\Omega\frac{|\nabla h_\varepsilon|^2}{h_\varepsilon}\right) \|\nabla \phi\|_{L^\infty(\Omega)}  \\[5pt]
  & \to 0 ~\text{as}~ \varepsilon = \varepsilon_j \searrow 0.
\end{align}
By \eqref{reg_assumptions} and \eqref{mth6}
\begin{equation}
  \label{hweak_sol1}
  - \int_0^T \int_\Omega h_{\varepsilon} \partial_t \phi - \int_\Omega h_{0\varepsilon} \phi(\cdot, 0) \to - \int_0^T \int_\Omega h \partial_t \phi - \int_\Omega h_{0} \phi(\cdot, 0),
\end{equation}
as $\varepsilon = \varepsilon_j \searrow 0$. Together,  \eqref{mth411} and \eqref{mth6} give us
\begin{equation}
\label{hweak_sol2}
  - \gamma_1\int_0^T \int_\Omega  h_\varepsilon c_{1\varepsilon} \phi \to -  \gamma_1 \int_0^T \int_\Omega h c_{1} \phi.
\end{equation}
as $\varepsilon = \varepsilon_j \searrow 0$. Taken together,  \eqref{mth4} and \eqref{mth6} result in
 \begin{equation*}
\label{hweak_sol3}
  -  \gamma_2 h_\varepsilon c_{2\varepsilon} \to -  \gamma_2 h c_{2}  \quad \text{in}~L^{\frac{5}{4}}(\Omega \times (0, T))
\end{equation*}
as $\varepsilon = \varepsilon_j \searrow 0$, which implies that
\begin{align}
 \label{hweak_sol4}
  - \gamma_2 \int_0^T \int_\Omega  h_\varepsilon c_{2\varepsilon} \phi \to -\gamma_2 \int_0^T \int_\Omega h c_{2} \phi.
\end{align}
 Taken together, \eqref{mth4} and \eqref{mth6} result in
 \begin{equation}
\label{hweak_sol5}
  -  \mu \int_0^T \int_\Omega  h_\varepsilon \phi + \int_0^T \int_\Omega \frac{c_{2\varepsilon}}{1 + c_{2\varepsilon}} \phi \to -  \mu \int_0^T \int_\Omega  h \phi + \int_0^T \int_\Omega \frac{c_{2}}{1 + c_{2}}\phi
\end{equation}
as $\varepsilon = \varepsilon_j \searrow 0$.
The convergence results \eqref{artificial1}-\eqref{hweak_sol5} enable us to pass the limit $\varepsilon = \varepsilon_j \searrow 0$ in \eqref{mthdefeq3} to have
\begin{align}
   \label{wmthdefeq22}
 - \int_0^T \int_\Omega h \partial_t \phi - \int_\Omega h_{0} \phi(\cdot, 0) & = - \gamma_1 \int_0^T \int_\Omega h c_{1} \phi  - \gamma_2 \int_0^T \int_\Omega h c_{2}\phi - \mu \int_0^T \int_\Omega  h \phi  +  \int_0^T \int_\Omega \dfrac{c_{2}}{1 + c_{2}} \phi.
\end{align}
Similarly, we can deduce
\begin{align}
   \label{mthdefeq4}
  - \int_0^T \int_\Omega \tau \partial_t \phi - \int_\Omega \tau_{0} \phi(\cdot, 0) = - \delta \int_0^T \int_\Omega \tau c_{1} \phi - \sigma \int_0^T \int_\Omega \tau \phi - \int_0^T \int_\Omega \dfrac{c_{2}}{1 + c_{2}} \phi.
\end{align}
By collecting \eqref{wmthdefeq1}, \eqref{mthdefeq2}, \eqref{wmthdefeq22}, and \eqref{mthdefeq4}, we complete the proof.

\end{proof}

\section*{Acknowledgement} 

\noindent
This work was funded by a DFG research grant within the SPP 2311 \textit{Robust coupling of continuum-biomechanical in silico models to establish active biological system models for later use in clinical applications - Co-design of modeling, numerics and usability} (grant nr. 679005).

%%%%%%%%%%%%%%%%%%%%%%%%%%%%%%%%%%%%%%%%%%%%%%%%%%%%%%%%%%%%%%%%%%%%%%%%%%%%%%%%%%%%%%%%%%%%%%%%%%%%%%%%%%%%%%%%%%%%%%%%%%%%%%%%%%%%%%%%%%%%%%%%%%%%%%%%%%%%%%%%%%%%%%%%%%%%%%%%%%%%%%%%%%%%%%%%%%%%%%%%%%%%%%%%%%%%

\phantomsection
\printbibliography

\end{document}